\documentclass[12pt]{article}
\usepackage{amssymb}
\usepackage{amsfonts}
\usepackage{mathrsfs}
\usepackage{amsmath}
\usepackage{pifont}
\usepackage{amsopn}
\usepackage{bm}
\usepackage{amscd}
\usepackage{framed}
\usepackage{delarray}
\usepackage{xypic}
\usepackage{exscale}
\usepackage{amsthm}
\usepackage{hyperref}
\input diagxy
\begin{document}
\newcommand{\Proof}{\noindent {\bf Proof.}}
\newcommand{\Section}[1]{
   \stepcounter{section}
   \bigskip\noindent
   {\bf\hbox{\thesection.~}#1}\par
   \nopagebreak
   \medskip
   \renewcommand{\theequation}{\thesection.\arabic{equation}}
   \setcounter{equation}{0}
   \setcounter{subsection}{0}}
\theoremstyle{definition}
\newtheorem{definition}{Definition}[section]
\newtheorem{remark}{Remark}[section]
\theoremstyle{plain}
\newtheorem{lemma}{Lemma}[section]
\newtheorem{theorem}{Theorem}[section]
\newtheorem{example}{Example}[section]
\newtheorem{corollary}{Corollary}[section]
\newtheorem{proposition}{Proposition}[section]
\newtheorem{question}{Question}[section]
\title{{\bf Partially Ordered Sheaves on a Locale. II}
\thanks{\it Project supported by NSFC (11171156) and PCSIRT (IRTO0742).}}
\author{{\small Wei He}\thanks{{{\it E-mail address:} weihe@njnu.edu.cn}}\\
 \small Institute
of Mathematics, Nanjing Normal University, Nanjing, 210097, China}

\date{}
\maketitle

\begin{abstract}
This paper is a continuation of our paper [1]. In this paper, we first introduce the concept of directed complete partially
ordered sheaves (shortly dcposheaves) on a given locale. some internal characterizations of dcposheaves and meet continuous dcposheaves are given respectively. We also give characterizations of continuous posheaves and  completely distributive posheaves,  and show that an algebraic completely distributive posheaf is spatial.

 \noindent{\bf Keywords:}  sheaf; partially ordered sheaf; continuous partially ordered sheaf.

\noindent{\bf Mathematics Subject Classifications}(2000):  18F20; 06A06.
\end{abstract}

\section{Introduction}

 In classical order algebraical theory, continuous dcpo (or domains) are very important class of posets. It also closely concerned with topology,
  i.e. a continuous frame is just a sober locally compact topology. In this
  paper we first introduce the concept of directed complete partially ordered sheaves
on a locale $X$.   Some internal characterizations of directed complete
partially ordered sheaves and meet continuous partially ordered
sheaves are presented respectively.  Then we introduce the concept
of continuous posheaves which are continuous directed complete
partially ordered objects in a localic topos $Sh(X)$.
   We show that continuous frame sheaves correspond to open fame homomorphisms satisfying an additional
   condition. We also show that an algebraic
completely distributive posheaf is spatial. Throughout this paper,
when we write $X$ for a locale, we will write $\mathcal{O}(X)$ for
the corresponding frame.

\section{Directed complete posheaves}

 Let $X$ be a locale and $F$ a  sheaf on $X$. Recall that a point of  $F$ to be a  morphism $p: \hat{1} \rightarrow F$ with $\hat{1}$ a subsheaf of
 the terminal sheaf $1$. A point of the form $1 \rightarrow F$  will be called  a global point of $F$. For a point $p: \hat{1} \rightarrow F$ of $F$,
 we write $dom(p)$ for the largest open $u \in \mathcal{O}(X)$ with $p(u) \not= \emptyset$, i.e. $dom(u) = \bigvee \{u \in \mathcal{O}(X)
 \mid p(u) \not= \emptyset\}$, and call it the domain of $p$. The set of all points of a sheaf $F$ will be denoted by $F_{p}$.

 Let $F$ be a posheaf on $X$, and $D \subseteq F_{p}$  a family of points. We call $D$  a directed family if for any two points $ p_{1}, p_{2} $ in $D$ with $dom(p_{1}) \leq dom(p_{2})$,
   there exists a point $p \in D$ such that $p_{1} \leq p$ and $p_{2} \leq p$.

 \begin{definition} A posheaf $F$ is said to be directed if there exists a directed family $D \subseteq F_{p}$ which can generate $F$. \end{definition}

 By the  above definition, we know that a posheaf $F$ is  directed if and only if for $\forall u \in \mathcal{O}(X)$, there exists
 a directed set $D(u) \subseteq F(u)$ such that for any two elements $x \in D(u),  y \in D(v)$ with $u \leq v$,  there exists an element $z \in D(w)$
 such that $x \leq z|_{u}, y \leq z|_{v} $,   and for $\forall x \in F(u)$ there is a cover $u = \bigvee u_{i}$ satisfying $x|_{u_{i}} \in D(u_{i})$.
 If a directed subsheaf $G$ of $F$
 is also a downsheaf then we call it an ideal of $F$.

  For a posheaf $F$ on a locale $X$, note that for any directed subsheaf $D$ of $F$, $D^{u}$ is still directed, hence we can construct
  a directed-powersheaf  $\mathcal{D}F$ of $F$ defined by $\mathcal{D}F(u) = D(F^{u})$ where  $D(F^{u})$ is the set of
   all directed subsheaves of $F^{u}$, and each restriction map
  $\mathcal{D}F(u) \rightarrow \mathcal{D}F(v)$ is same as the restriction map for the power sheaf
  $\mathbb{P}F$. So $\mathcal{D}F$ is a subsheaf of $\mathbb{P}F$.

 If $p$ is a point of a posheaf $F$, the principle ideal $\downarrow p$ is clearly a directed subsheaf of $F$. Hence the principle ideal embedding
 morphism $\downarrow: F \rightarrow \mathbb{P}F$ can be factored through  $\mathcal{D}F$,
    \[\downarrow: F \rightarrow \mathcal{D}F\]

\begin{definition}  A posheaf $F$ is called directed complete if the principle ideal embedding $\downarrow: F \rightarrow \mathcal{D}(F)$ has a left adjoint.\end{definition}

\begin{proposition} For any posheaf $F$ on a locale $X$, the following are equivalent:

(1) $F$ is directed complete.

(2)  For every directed subsheaf $D$ of $F$, $\bigvee D$ exists  and
it can be extended to a global point $p$ of $F$ such that for each
$u \in \mathcal{O}(X)$, $p|_{u}$ is the least element in $F(u)$
satisfying $D^{u} \subseteq \downarrow p|_{u}$.

(3)  For every ideal $D$ of $F$, $\bigvee D$ exists  and it can be
extended to a global point $p$ of $F$ such that for each $u \in
\mathcal{O}(X)$, $p|_{u}$ is the least element in $F(u)$ satisfying
$D^{u} \subseteq \downarrow p|_{u}$.

(4) Every restriction map is surjective. For every $u \in
\mathcal{O}(X), F(u)$ is a directed complete poset. Every
restriction map $F(v) \rightarrow F(u)$ preserves directed joins and
has a left adjoint.
\end{proposition}

{\bf Proof}\, $(1) \Leftrightarrow (2)$ is same as the proof of
Proposition 3.2 in [1].

 $(2) \Rightarrow (3)$ is clear.

 $(3) \Rightarrow (2)$ Suppose $D$ is a directed subsheaf of $F$,
  then $\downarrow D$ is also directed hence an ideal and $\bigvee D = \bigvee \downarrow
  D$. By (3) $\bigvee D$ can be
extended to a global point $p$ of $F$ such that for each $u \in
\mathcal{O}(X)$, $p|_{u}$ is the least element in $F(u)$ satisfying
$(\downarrow D)^{u} \subseteq \downarrow p|_{u}$. We note that
$(\downarrow D)^{u} = \downarrow D^{u}$, hence (2) holds.

 $(1) \Rightarrow (4)$  To show each restriction map $F(1_X) \rightarrow F(u)$ is surjective, we only need to note that
  for any $x \in F(u)$, regarded as a point of $F$, it is directed. Suppose $D \subseteq F(u)$ is directed.
    Denote $\bar{D}$ the presheaf generated by $D$: \[\bar{D}(v) = \left\{\begin{array}{cc}
\{x \in F(v)\mid \exists y \in D,  x = y|_{v}\}, & v \leq u\\
\emptyset, & v \not\leq u
\end{array}
\right.\]  Write $\hat{D}$ the sheaf generated by $\bar{D}$, then $
\hat{D}$ is directed and  $\bigvee \hat{D}$ exists. Similar to the
proof of proposition 3.2 in [1],  we can show that $\bigvee \hat{D}$ is
also the least upper bound of $D$ in the poset $F(u)$,  and for $v \leq u$,  we have $(\bigvee
D)|_{v} = \bigvee D|_{v}$. Hence $F(u)$ is directed complete and the
restriction map $F(u) \rightarrow F(v)$ preserves directed joins.  Moreover,  $F(u) \rightarrow F(v)$ has a left adjoint.

$(4) \Rightarrow (2)$  Suppose $A$ is a directed subsheaf of $F$,
and $D \subseteq F_{p}$ is a directed family which can generate $A$.
Regard $D$ as a sub-presheaf of $F$, we know $D$ and $A$ have the
same upper bound set. For every $u \in \mathcal{O}(X)$, write $D_{u}
= \{x \in F(u) \mid \exists p \in D, x = p|_{u}\}$. Then $D_{u}$ is
directed in $F(u)$ hence $\bigvee D_{u}$ exists.  For $v \leq u$ in
$\mathcal{O}(X)$, $(\bigvee D_{u})|_{v} = \bigvee D_{u}|_{v} \leq
\bigvee D_{v}$ since the restriction map $F(u) \rightarrow F(v)$
preserves directed joins.  For any $x \in D_{v}$, by the directness
of $D$ there exists $y \in D_{u}$ such that $x \leq y|_{v}$, hence
the converse inequality $\bigvee D_{v} \leq (\bigvee D_{u})|_{v} $
holds. This shows that $\{\bigvee D_{u} | D_{u} \not= \emptyset\}$
is a compatible family, hence it can be patched to an unique element
$p \in F(\bigvee \{u \in \mathcal{O}(X) \mid D_{u}
\not=\emptyset\})$. It is clear that $p$ is the least upper bound of
$A$. The point $p$  can be extended to a global point $\bar{p}$ of
$F$ such that for each $u \in \mathcal{O}(X)$, $\bar{p}|_{u}$ is the
least element in $F(u)$ satisfying $D^{u} \subseteq \downarrow
\bar{p}|_{u}$ since the restriction map $F(1_{X}) \rightarrow
F(\bigvee \{u \in \mathcal{O}(X) \mid D_{u} \not=\emptyset\})$ has a
left adjoint. $\hfill\Box$

\begin{corollary} Let $X$ be a locale and $F$ a posheaf on $X$. The
following conditions are equivalent:

(1) $F$ is complete.

(2) $F$ is finite upper-complete and directed complete.
\end{corollary}

 We write
$DCPO^{\ast}$  for the category of all directed complete partially
ordered sets and surjective maps having left adjoint and preserving
directed joins.  Then we have

\begin{corollary} Let $X$ be a locale. $F$ is a directed complete posheaf on $X$ if and only if $F$ is a sheaf over $DCPO^{\ast}$ and satisfying the condition
(POS3).\end{corollary}

Now we consider those morphisms which preserve directed sups between
directed complete posheaves.  suppose $\alpha: F \rightarrow G$ is
an order-preserving morphism and $S$ is a directed subsheaf of $F$.
If $D \subseteq F_{p}$ is a directed family which generates $S$,
then the image $\alpha D$ of $D$ under $\alpha$ is a directed family
which generates the image $\hat{\alpha S}$ of $S$, hence
$\hat{\alpha S}$ is directed.  Let $F$ and $G$ be two directed
complete posheaves on a locale $X$ and $\alpha: F \rightarrow G$ be
an order-preserving morphism. Sine the image of a directed subsheaf
is directed, we can construct an order-preserving morphism
$\alpha_{\ast}: \mathcal{D}F \rightarrow \mathcal{D}G$ such that for
each $u \in \mathcal{O}(X)$ and  $S \in D(F^{u})$,
${\alpha_{\ast}}_{u}: D(F^{u}) \rightarrow D(G^{u})$ maps $S$ to the
image of $S$ under $\alpha$.

\begin{definition} Let $F$ and $G$ be two directed
complete posheaves on a locale $X$ and $\alpha: F \rightarrow G$ be
an order-preserving morphism.  We call $\alpha$ preserving directed
sups if the following square commutes:
 \[
\xymatrix{\mathcal{D}F \ar[r]^{\alpha_{\ast}}\ar[d]_{sup_{F}}&
\mathcal{D}G
\ar[d]^{sup_{G}}\\
      F  \ar[r]^{\alpha}& G}
\] where $sup_{F}: \mathcal{D}(F) \rightarrow F$ and $sup_{G}: \mathcal{D}(G) \rightarrow
G$ are respectively the left adjoint of the principle ideal embedding
$F \rightarrow \mathcal{D}F$ and $G \rightarrow \mathcal{D}G$.
\end{definition}

Similar to the proof of proposition 3.3 in [1], we have the following
result.

\begin{proposition}   Let $F$ and $G$ be directed complete posheaves on a locale $X$ and let $\alpha: F
\rightarrow G$ be an order-preserving morphism. The following
conditions are equivalent:

(1)  $\alpha: F \rightarrow G$ preserves directed sups.

(2) For each $u \in \mathcal{O}(X)$, $\alpha_{u}: F(u) \rightarrow
G(u)$ preserves directed joins and the following square commutes for
any $v \leq u$ in $\mathcal{O}(X)$
 \[
\xymatrix{ F(u)\ar[r]^{\alpha_{u}}&  G(u)\\
           F(v)\ar[r]^{\alpha_{v}}\ar[u]^{f_{uv}}
           & G(v)\ar[u]^{g_{uv}}}
\] where $f_{uv}: F(v) \rightarrow F(u)$ and $g_{uv}: G(v) \rightarrow G(u)$ are the left adjoint of the restriction maps $F(u) \rightarrow
F(v)$ and $G(u) \rightarrow G(v)$ respectively.
\end{proposition}

Now we consider the generalization of meet continuous semilattices.  Recall that a  posheaf $F$ is said to be finite inf-complete if  $F \rightarrow 1 $ and the diagonal $F \rightarrow F\times F$ both has a right adjoint. An finite inf-complete posheaf $F$ cab be characterized as a posheaf $F$ such that for every  $u \in \mathcal{O}(X)$, $F(u)$ is a inf-semilattice and every restriction map $F(u) \rightarrow F(v)$ preserves finite meets.
Let $F$ be a directed complete finite inf-complete posheaf.  We define a meet morphism $\mu_{F}: F \times
\mathbb{P}F \rightarrow \mathbb{P}F$ of points with subsheaves such that for
each $u \in \mathcal{O}(X)$, $x \in F(u)$, and $S \in sub(F^{u})$,
$\mu_{F}(u)(x, S)$ be the directed subsheaf of $F^{u}$ generated by the
sub-presheaf $\bar{S}$ of $F^{u}$ with $\bar{S}(v) = \{x|_{v} \wedge
y \mid y \in S(v)\}$ for each $v \leq u$. Note that since each
restriction map preserves meets so $\bar{S}$ is indeed a
sub-presheaf of $F^{u}$.

\begin{definition} Let $F$ be a  directed complete finite inf-complete posheaf on a locale $X$.
$F$ is  said to be a meet continuous posheaf  if the
following square commutes:

 \[
\xymatrix{F \times \mathbb{P}F \ar[r]^{1_{F}\times
sup_{F}}\ar[d]_{\mu_{F}}&  F \times F
\ar[d]^{m_{F}}\\
      \mathbb{P}F   \ar[r]^{sup_{F}}& F}
\]
 where $m_{F}: F \times F \rightarrow F$ is the
right adjoint of the diagonal $F \rightarrow F \times F$.
\end{definition}

Similar to the proof of Proposition 3.7 in [1], we have the following result.

\begin{proposition}  Let $F$ be a directed complete finite inf-complete posheaf on a locale $X$.
The following conditions are equivalent:

(1) $F$ is a meet continuous posheaf .

(2) Each $F(u)$ is a meet continuous semilattice for $u \in
\mathcal{O}(X)$, and if $v \leq u$ in $\mathcal{O}(X)$,  then $x
\wedge l_{uv}(y) = l_{uv}( x|_{v} \wedge y)$ holds for any $x \in
F(u), y \in F(v)$ where $l_{uv}: F(v) \rightarrow F(u)$ is the left
adjoint of the restriction map $F(u) \rightarrow F(v)$.
\end{proposition}

\begin{corollary} Let $F$  be a complete  posheaf on a locale $X$. The following conditions are equivalent:

(1) $F$ is a complete Heyting sheaf.

(2) $F$ is a meet continuous distributive posheaf.
\end{corollary}

 \section{Continuous Posheaves}

  For a posheaf $F$ on a locale $X$,  we can construct
  an ideal-powersheaf  $\mathcal{ID}F$ of $F$ defined by $\mathcal{ID}F(u) = ID(F^{u})$ where  $ID(F^{u})$ is the set of
   all ideals of $F^{u}$, and each restriction map
  $\mathcal{ID}F(u) \rightarrow \mathcal{ID}F(v)$ is same as the restriction map for
  $\mathbb{P}F$. So $\mathcal{ID}F$ is a subsheaf of $\mathbb{P}F$.

 If $p$ is a point of a posheaf $F$, the principle ideal $\downarrow p$ is clearly an ideal of $F$. Hence the principle ideal embedding
 morphism $\downarrow: F \rightarrow \mathbb{P}F$ can be factored through  $\mathcal{ID}F$,
    \[\downarrow: F \rightarrow \mathcal{ID}F\]

By Proposition 2.1,  we have the following result.

\begin{lemma}  For any posheaf $F$ on a locale $X$,  $F$ is directed complete if an only if the principle ideal embedding $\downarrow: F \rightarrow \mathcal{ID}(F)$ has a left adjoint.
\end{lemma}

Let $F$ be a directed complete posheaf on $X$  and $sup_{F}: \mathcal{ID}(F) \rightarrow F$  the left adjoint of the principle ideal embedding $\downarrow: F \rightarrow \mathcal{ID}(F)$. $F$ is said to be a continuous posheaf if $sup_{F}: \mathcal{ID}(F) \rightarrow F$ has a left adjoint.

\begin{proposition} Let $F$ be a directed complete posheaf on $X$. The follows are equivalent:

(1) F is a continuous posheaf;

(2) for each global point $t \in F(1_{X})$,  there exists a smallest ideal $I_{t}$ of $F$ such that $\bigvee I_{t} = t$ and for each $u \in \mathcal{O}(X)$,  the restriction $I_{t}^{u}$ of $I_{t}$ to $F^{u}$ be the smallest ideal such that
$\bigvee I_{t}^{u} = t|_{u}$;

(3) $sup_{F}: \mathcal{ID}(F) \rightarrow F$  is an inf-preserving morphism.  \end{proposition}

\begin{corollary} Let $F$  be a continuous posheaf on $X$. Then for each $y \in F(u)$,  $u \in \mathcal{O}(X)$, there exists $min\{v \leq u \mid l_{uv}(y|_{v}) = y\}$ which we denoted by $u(y)$ with the property that for each $v \leq u$, $v(y|_{v}) = v \wedge u(y)$ where $l_{uv}: F(v) \rightarrow F(u)$ be the left adjoint of the restriction map $F(u) \rightarrow F(v)$. \end{corollary}

\begin{proposition} Let $F$ be a continuous posheaf and $G$ a directed complete posheaf on $X$.  If $\alpha: F \rightarrow G$ is a surjective morphism which preserves directed sups and has a left adjoint then $G$ is continuous. \end{proposition}

{\bf Proof}\, Write $\beta: G \rightarrow F$ to be the left adjoint of $\alpha$ and $\gamma: F \rightarrow \mathcal{ID}(F)$ to be the left adjoint of the morphism $sup_{F}: \mathcal{ID}(F) \rightarrow F$. We now show that the composite $\alpha_{\ast} \gamma \beta: G \rightarrow F \rightarrow \mathcal{ID}(F) \rightarrow \mathcal{ID}(G)$ is the left adjoint of the morphism
$sup_{G}: \mathcal{ID}(G) \rightarrow G$ where $\alpha_{\ast}:  \mathcal{ID}(F) \rightarrow \mathcal{ID}(G)$ be the image morphism determined by $\alpha$.

Suppose $u \in \mathcal{O}(X)$, and $y \in G(u)$. If $\alpha_{\ast u}(\gamma_{u}(\beta_{u}(y))) \subseteq J$ for some $J \in \mathcal{ID}(G)^{u}$, then $sup_{G u}(\alpha_{\ast u}(\gamma_{u}(\beta_{u}(y)))) \leq sup_{G u}(J)$. But $sup_{G u}(\alpha_{\ast u}(\gamma_{u}(\beta_{u}(y)))) = \alpha_{u}(sup_{F u}(\gamma_{u}(\beta_{u}(y)))) = \alpha_{u}(\beta\beta_{u}(y)) = y$ since $\alpha$ preserves directed sups and is a surjective morphism.

Conversely, if $y \leq sup_{G u}(J)$ for some $J \in \mathcal{ID}(G)^{u}$, then $\beta_{u}(y) \leq \beta_{u}(sup_{G u}(J)) = sup_{G u}(\beta_{\ast u}(J))$ since $\beta$ preserves sups. Hence $\gamma_{u}(\beta_{u}(y)) \subseteq \beta_{\ast u}(J)$ by the adjointness. It implies that $\alpha_{\ast u}(\gamma_{u}(\beta_{u}(y))) \subseteq \alpha_{\ast u}(\beta_{\ast u}(J)) = J$.  $\hfill\Box$

\begin{definition} Let $F$ be a complete posheaf.  If $sup_{F}: \mathbb{D}(F) \rightarrow F$ has a left adjoint then we call $F$ a completely distributive posheaf, where $sup_{F}: \mathbb{D}(F) \rightarrow F$ is the left adjoint of the principle ideal embedding
$F \rightarrow \mathbb{D}F$. \end{definition}

\begin{proposition} Let $F$ be a complete posheaf on $X$. The follows are equivalent:

(1) F is a completely distributive posheaf;

(2) For each global point $t \in F(1_{X})$,  there exists a lowersheaf $G_{t}$ of $F$ such that $\bigvee G_{t} = t$ and for each $u \in \mathcal{O}(X)$,  the restriction $G_{t}^{u}$ of $G_{t}$ to $F^{u}$ be the smallest lowersheaf such that
$\bigvee G_{t}^{u} = t|_{u}$;

(3) $sup_{F}: \mathbb{D}(F) \rightarrow F$  is an inf-preserving morphism.  \end{proposition}

\begin{lemma} Every completely distributive posheaf is distributive. \end{lemma}

{\bf Proof}\, Let  $F$ be a completely distributive posheaf on $X$. We only need to show that the global sections set $F(1_{\mathcal{O}(X)})$ is distributive.

Suppose $a, b ,c \in F(1_{\mathcal{O}(X)})$, write $t = a \wedge(b \vee c)$. Denote $G_{t}$ the smallest lowersheaf of $F$ such that $\bigvee G_{t} = t$ and $\widetilde{bc}$ be the lowersheaf of $F$ generated by $\{b, c\}$. Then we have $G_{t} \subseteq \widetilde{bc}$ by the adjointnees and the fact that $t \leq b \vee c = \bigvee \widetilde{bc}$. Hence  $G_{t} \subseteq \widetilde{(a\wedge b)(a\wedge c)}$ since $G_{t} \subseteq \downarrow a$, where $\widetilde{(a\wedge b)(a\wedge c)}$ be the lowersheaf of $F$ generated by $\{a \wedge b, a \wedge c\}$. This shows that $t = \bigvee G_{t} \leq  \bigvee \widetilde{(a\wedge b)(a\wedge c)} = (a\wedge b) \vee (a\wedge c)$. Thus $a \wedge(b \vee c) =  (a\wedge b) \vee (a\wedge c)$.  $\hfill\Box$

Now we construct a left adjoint for the inclusion $\mathcal{ID}F \rightarrow \mathbb{D}F$.

For each $u \in \mathcal{O}(X)$ and a subsheaf $G$ of $F^{u}$, consider the sub-presheaf $\bar{G}$ of $F^{u}$ defined by $\bar{G}(v) = \{l_{w_{1}v}(x_{1})\wedge \cdot\cdot\cdot \wedge l_{w_{n}v}(x_{n}) \mid x_{i} \in G(w_{i}), w_{i} \leq v, i= 1, \cdot\cdot\cdot, n\}$ where $l_{w_{i}v}: G(w_{i}) \rightarrow G(v)$ be the left adjoint of the restriction map $G(v) \rightarrow G(w_{i}))$. Let $\hat{G}$ be the ideal generated by $\bar{G}$. By Lemma 3.1 in [1], It is straightforward to show that the correspondence $ \mathcal{D}F(u) \rightarrow \mathcal{ID}F, G \mapsto \hat{G}$ determines a functor which is left adjoint to the  the inclusion $\mathcal{ID}F \rightarrow \mathcal{D}F$.

 We have the following diagram commutes
\[
\xymatrix{
  \mathcal{ID}F \ar[rr]^{}\ar[dr]_<<<<<{sup}& & \mathcal{D}F \ar[dl]^<<<<<{sup}\\
 &  F  &}
\]
By the above arguments we have the following result.

\begin{lemma} Every completely distributive posheaf is a continuous posheaf. \end{lemma}

\section{Continuous Frame Sheaves}

\begin{proposition} Let $X$ be a locale and and let $F$ be a complete posheaf on $X$. If $F$ is continuous and distributive and for $v \leq u$ in $\mathcal{O}(X)$,  we have $x
\wedge l_{uv}(y) = l_{uv}( x|_{v} \wedge y)$ holds for any $x \in
F(u), y \in F(v)$ where $l_{uv}: F(v) \rightarrow F(u)$ is the left
adjoint of the restriction map $F(u) \rightarrow F(v)$. Then $F$ is a frame sheaf. \end{proposition}

{\bf Proof}\, We only need to show that each $F(u)$ is a frame for $u \in \mathcal{O}(X)$. Let $x \in F(u)$ and $D \subset F(u)$ be a directed set. Write $t = x \wedge \bigvee D$,  denote $I_{t}$  the smallest ideal on $F^{u}$ such that $\bigvee I_{t} = t$ and $\bar{D}$  the ideal on $F^{u}$ generated by $D$. Then we have $I_{t} \subset \bar{D}$ by the adjointness. Hence for each $v \leq u$ and $y \in I_{t}(v)$ there exists a cover $v = \bigvee v_{j}$ such that $y|_{v_{j}} \leq d_{j}|_{v_{j}}$ for some $d_{j} \in D$. This implies that $y \leq x \wedge \bar{d}$ for some $\bar{d} \in \bar{D}$. Thus $t \leq \bigvee xD$. Note that the distributivity of $F$ implies that each $F(u)$ is a distributive lattice, so $F(u)$ is a frame.  $\hfill\Box$

Let $f: Y \rightarrow X$ be a locale morphism. By the equivalence of the category of frame sheaves on $X$ and the  category of frames under $\mathcal{O}(X)$, we can regard $f: Y \rightarrow X$ as a frame sheaf. Now it is natural to ask under what conditions $f: Y \rightarrow X$ corresponds to continuous frame sheaf.

\begin{definition} Let $f: Y \rightarrow X$ be a locale morphism such that the corresponding frame homomorphism $f^{\ast}: \mathcal{O}(X) \rightarrow \mathcal{O}(Y)$
has a left adjoint $f^{\sharp}: \mathcal{O}(Y) \rightarrow \mathcal{O}(X)$.  $y, z \in \mathcal{O}(Y)$. We call that $z$ is $f-way bellow$ $y$, written as  $z \ll_{f} y$, if and only if $z \leq f^{\ast}f^{\sharp}(y)$,  and for all $u \in \mathcal{O}(X)$, any directed subset $D \subset \mathcal{O}(Y)$ with $y \wedge f^{\ast}(u) \leq \bigvee D$,  there exists a cover $u \wedge f^{\sharp}(z) = \bigvee \{u_{j} \in \mathcal{O}(X)\mid j \in J\}$ such that for $\forall j \in J$, $z \wedge f^{\ast}(u_{j}) \leq d_{j} \wedge f^{\ast}(u_{j})$ for some $d_{j} \in D$. \end{definition}

For a given locale $f: Y \rightarrow X$ over $X$, the above binary relation $\ll_{f} $ is a generalization of the way-below relation $\ll$ on $ \mathcal{O}(Y)$. Indeed, If we take $X$ to be the terminal locale, i.e. $\mathcal{O}(X) = 2$ be the initial frame, then $z \ll_{f} y$  if and only if $z \ll y$. Similar to the way-below relation, $\ll_{f}$ relation has following properties.

\begin{proposition} Let $f: Y \rightarrow X$ be a locale morphism such that the corresponding frame homomorphism $f^{\ast}: \mathcal{O}(X) \rightarrow \mathcal{O}(Y)$
has a left adjoint $f^{\sharp}: \mathcal{O}(Y) \rightarrow \mathcal{O}(X)$.   The following statement hold for all $x, y, z \in \mathcal{O}(Y)$:

(1) $z \ll_{f} y$ implies $z \leq y$;

(2) $x \leq z \ll_{f} y \leq t$ implies $x \ll_{f} t$;

(3) $z \ll_{f} y$ and $x \ll_{f} y$ implies $x \vee z \ll_{f} y$;

(4) $0 \ll_{f} x$;

(5) If $f: Y \rightarrow X$ is an open morphism then $z \ll_{f} y$ implies $z \wedge f^{\ast}(u) \ll_{f} y \wedge f^{\ast}(u)$ for all $u \in \mathcal{O}(X)$.

(6)  Suppose $z \leq f^{\ast}(u)$ and $u = \bigvee_{i\in I} u_{i}$. If $z \wedge f^{\ast}(u_{i}) \ll_{f} y_{i}$ for each $i \in  I$ then $z \ll_{f} \bigvee_{i\in I} y_{i}$.
 \end{proposition}

 {\bf Proof}\, (1) Take $u = f^{\sharp}(y)$ then $y \wedge f^{\ast}(u) = y$. We have $f^{\sharp}(z) \leq u$ since $z \leq f^{\ast}f^{\sharp}(y) = f^{\ast}(u)$.  Let $D = \{y\}$ then there exists an cover $f^{\sharp}(z) = \bigvee \{u_{j} \in \mathcal{O}(X)\mid j \in J\} $ such that for $\forall j \in J$,  $z \wedge f^{\ast}(u_{j}) \leq y  \wedge f^{\ast}(u_{j})$. Hence
 $z \wedge \bigvee_{j\in J} f^{\ast}(u_{j}) \leq y  \wedge \bigvee_{j\in J} f^{\ast}(u_{j})$. This implies  that $z = z \wedge f^{\ast}f^{\sharp}(z) = z \wedge f^{\ast}(\bigvee_{j\in J} u_{j}) \leq y \wedge f^{\ast}(\bigvee_{j\in J} (u_{j}))  \leq y$.

 (2) Clear.

 (3) Let $u \in \mathcal{O}(X)$ and  $D \subset \mathcal{O}(Y)$ be a directed subset such that $y \wedge f^{\ast}(u) \leq \bigvee D$. There exists a cover $u \wedge f^{\sharp}(z) = \bigvee \{u_{j} \in \mathcal{O}(X)\mid j \in J\} $ such that for $\forall j \in J$, $z \wedge f^{\ast}(u_{j}) \leq d_{j} \wedge f^{\ast}(u_{j})$ for some $d_{j} \in D$ and another cover $u \wedge f^{\sharp}(x) = \bigvee \{v_{i} \in \mathcal{O}(X)\mid i \in I\} $ such that for $\forall i \in I$, $x \wedge f^{\ast}(v_{i}) \leq d_{i} \wedge f^{\ast}(v_{i})$ for some $d_{i} \in D$. Then $u \wedge f^{\sharp}(x \vee z) = \bigvee \{v_{i} \vee u_{j} \mid i \in I, j\in J\}$, and $(x \vee z) \wedge f^{\ast}(v_{i} \vee u_{j}) \leq ((z \vee d_{i}) \wedge f^{\ast}(v_{i})) \vee ((x \vee d_{j}) \wedge f^{\ast}(u_{j})) \leq (z \vee d_{i} \vee x \vee d_{j}) \wedge f^{\ast}(u_{i} \vee v_{j}) \leq (x \vee z \vee d_{k}) \wedge  f^{\ast}(u_{i} \vee v_{j})$ for some $d_{k} \in D$ since $D$ is directed.

 (4) Clear.

 (5) Let $u \in \mathcal{O}(X)$. First $z \leq f^{\ast}f^{\sharp}(y)$ implies that $z \wedge f^{\ast}(u) \leq  f^{\ast}f^{\sharp}(y) \wedge f^{\ast}(u) = f^{\ast}f^{\sharp}(y \wedge f^{\ast}(u))$.
Suppose  $y \wedge f^{\ast}(u) \wedge f^{\ast}(v)  \leq \bigvee D$ for a directed subset $D \subset \mathcal{O}(Y)$ and $v \in \mathcal{O}(X)$.  There exists a cover $u \wedge v \wedge f^{\sharp}(z) = \bigvee \{v_{j} \in \mathcal{O}(X)\mid j \in J\} $ such that for $\forall j \in J$, $z \wedge f^{\ast}(v_{j}) \leq d_{j} \wedge f^{\ast}(v_{j})$ for some $d_{j} \in D$. Thus $z \wedge f^{\ast}(u) \wedge f^{\ast}(v_{j}) \leq d_{j} \wedge f^{\ast}(v_{j})$ for all $j \in J$.

(6) Write $y = \bigvee y_{i}$ and $z_{i} = z \wedge f^{\ast}(u_{i})$ for $i \in I$. Suppose $y \wedge f^{\ast}(v) \leq  \bigvee D$ for some $v \in \mathcal{O}(X)$ and a directed subset $D \subset \mathcal{O}(Y)$, Then for each $i \in I$, there exists a cover $f^{\sharp}(z_{i}) \wedge v = \bigvee_{j\in J_{i}}  v^{i}_{j}$,   such that for each $j \in J_{i}$,  $z_{i} \wedge f^{\ast}(v^{i}_{j}) \leq d_{j} \wedge  f^{\ast}(v^{i}_{j})$ for some $d_{j} \in D$. It is clear that the cover  $f^{\sharp}(z) \wedge v = \bigvee_{i\in I} \bigvee_{j\in J_{i}}  v^{i}_{j}$  satisfies the condition.     $\hfill\Box$

 Write $\Downarrow_{f} y = \{z \in \mathcal{O}(Y) \mid z \ll_{f} y\}$. By Proposition 3.2 (6), we have the following result.

 \begin{corollary} Let $f: Y \rightarrow X$ be a locale morphism. The assignment $u \mapsto \Downarrow_{f} f^{\ast}(u)$, $u \in \mathcal{O}(X)$  defines a subsheaf of the frame sheaf  $f: Y \rightarrow X$. \end{corollary}

\begin{theorem} Let $f: Y \rightarrow X$ be a locale morphism. The corresponding frame sheaf of $f$ is a continuous frame sheaf if and only if $f: Y \rightarrow X$ is an open morphism and $y = \bigvee \{x \in \mathcal{O}(Y) \mid x  \ll_{f} y\}$ for all $y \in \mathcal{O}(Y)$. \end{theorem}

{\bf Proof}\, Suppose $f$  corresponds to a continuous frame sheaf $F$. For each $y \in \mathcal{O}(Y)$,  by corollary 3.1, there exists $min\{v \in \mathcal{O}(X)\mid y \leq f^{\ast}(v)\}$, hence $f^{\ast}: \mathcal{O}(X) \rightarrow \mathcal{O}(Y)$ has a left adjoint $f^{\sharp}$. Also $f^{\sharp}(y \wedge f^{\ast}(u)) = f^{\sharp}(y) \wedge u$ for all $y \in \mathcal{O}(Y)$ and $u \in \mathcal{O}(X)$ by corollary 3.1.

For each $y \in \mathcal{O}(Y)$, there exists a  smallest ideal $I(y)$ of $F$ such that $\bigvee I(y) = y$ and for each $u \in \mathcal{O}(X)$,  the restriction $I(y)^{u}$ of $I(y)$ to $F^{u}$ be the smallest ideal such that
$\bigvee I(y)^{u} = y \wedge f^{\ast}(u)$.  Now we show that $z \ll_{f} y$ for each $z \in I(y)_{p}$. Suppose that $z \in I(y)_{p}$, it is clear that $z \leq f^{\ast}f^{\sharp}(y)$.  Let $u \in \mathcal{O}(X)$ and $y \wedge f^{\ast}(u) \leq \bigvee D$ for a directed subset $D \subset \mathcal{O}(Y)$. Write $\bar{D}$ for the ideal generated by $D$.  Then we have $y \wedge f^{\ast}(u) \leq \bigvee \bar{D}$. Hence $I(y)^{u} \subset \bar{D}$.  It implies that $z \wedge f^{\ast}(u) \in \bar{D}$. Thus there exists a cover $u \wedge f^{\sharp}(z) = \bigvee \{u_{j} \in \mathcal{O}(X)\mid j \in J\} $ such that for $\forall j \in J$, $z \wedge f^{\ast}(u_{j}) \leq d_{j} \wedge f^{\ast}(u_{j})$ for some $d_{j} \in D$.

Conversely, suppose $f: Y \rightarrow X$ is an open morphism and $y = \bigvee \{x \in \mathcal{O}(Y) \mid x  \ll_{f} y\}$ for all $y \in \mathcal{O}(Y)$. For $y \in \mathcal{O}(Y)$, by proposition 4.2, the set $\{x \in \mathcal{O}(Y) \mid x  \ll_{f} y\}$ is directed. Write $I(y)$ for the ideal generated by the directed set $\{x \in \mathcal{O}(Y) \mid x  \ll_{f} y\}$. Then it is not difficult to show that $I(y)$ is the smallest ideal of $f: Y \rightarrow X$ such that $\bigvee I(y) = y$ and for each $u \in \mathcal{O}(X)$,  the restriction $I(y)^{u}$ of $I(y)$  be the smallest ideal such that
$\bigvee I(y)^{u} = y \wedge f^{\ast}(u)$.  $\hfill\Box$

Now we consider the conditions that makes a frame sheaf to be a completely distributive frame sheaf.

Let $f: Y \rightarrow X$ be a locale morphism such that the corresponding frame homomorphism $f^{\ast}: \mathcal{O}(X) \rightarrow \mathcal{O}(Y)$
has a left adjoint $f^{\sharp}: \mathcal{O}(Y) \rightarrow \mathcal{O}(X)$.  $y, z \in \mathcal{O}(Y)$. We introduce a stronger binary relation on $\mathcal{O}(Y)$ as following.

 We write  $z \lhd_{f} y$, if and only if $z \leq f^{\ast}f^{\sharp}(y)$,  and for all $u \in \mathcal{O}(X)$, any subset $B \subset \mathcal{O}(Y)$ with $y \wedge f^{\ast}(u) \leq \bigvee B$,  there exists a cover $u \wedge f^{\sharp}(z) = \bigvee \{u_{j} \in \mathcal{O}(X)\mid j \in J\} $ such that for $\forall j \in J$, $z \wedge f^{\ast}(u_{j}) \leq b_{j} \wedge f^{\ast}(u_{j})$ for some $b_{j} \in B$. Similar to the proof of proposition 4.2, we have the following result.

\begin{lemma} Let $f: Y \rightarrow X$ be a locale morphism such that the corresponding frame homomorphism $f^{\ast}: \mathcal{O}(X) \rightarrow \mathcal{O}(Y)$
has a left adjoint $f^{\sharp}: \mathcal{O}(Y) \rightarrow \mathcal{O}(X)$.   Then the following statements hold for all $x, y, z \in \mathcal{O}(Y)$:

(1) $z \lhd_{f} y$ implies $z \ll_{f} y$;

(2) $x \leq z \lhd_{f} y \leq t$ implies $x \lhd_{f} t$;

(3) $0 \lhd_{f} x$;

(4) If $f: Y \rightarrow X$ is an open morphism then $z \lhd_{f} y$ implies $z \wedge f^{\ast}(u) \lhd_{f} y \wedge f^{\ast}(u)$ for all $u \in \mathcal{O}(X)$.

(5)  Suppose $z \leq f^{\ast}(u)$ and $u = \bigvee_{i\in I} u_{i}$. If $z \wedge f^{\ast}(u_{i}) \lhd_{f} y_{i}$ for each $i \in  I$ then $z \lhd_{f} \bigvee_{i\in I} y_{i}$.

 \end{lemma}

Similar to the proof of theorem 4.1, the following result is clear.

\begin{theorem} Let $f: Y \rightarrow X$ be a locale morphism. The corresponding frame sheaf of $f$ is a completely distributive posheaf if and only if $f: Y \rightarrow X$ is an open morphism and $y = \bigvee \{x \in \mathcal{O}(Y) \mid x  \lhd_{f} y\}$ for all $y \in \mathcal{O}(Y)$. \end{theorem}

\begin{lemma} If $f: Y \rightarrow X$ is an open inclusion then $y \lhd_{f} y$ for all $y \in \mathcal{O}(Y)$. \end{lemma}

{\bf Proof}\,  Suppose $f: Y \rightarrow X$ be an open inclusion then $f^{\ast}:  \mathcal{O}(X) \rightarrow  \mathcal{O}(Y)$ is surjective. Hence for each $u \in  \mathcal{O}(X)$ and any subset $B \subset \mathcal{O}(Y)$ with $y \wedge f^{\ast}(u) \leq \bigvee B$, we have $u \wedge f^{\sharp}(y) = f^{\sharp}(y \wedge f^{\ast}(u)) \leq \bigvee f^{\sharp}(B)$.  This cover satisfying $y \wedge f^{\ast}f^{\sharp}(b) = y \wedge b \leq b$ for all $b \in B$.  $\hfill\Box$

\begin{proposition} If $f: Y \rightarrow X$ is a local homeomorphism then it corresponds to a completely distributive posheaf, hence a continuous posheaf. \end{proposition}

{\bf Proof}\,  We first note that for each open inclusion $i_{x}: \downarrow x \rightarrow Y$, $z \lhd_{fi_{x}} y$ implies that $z \lhd_{f} y$ for all $z, y \leq x$
since $(fi_{x})^{\sharp}(y) = f^{\sharp}(y)$ and $y \wedge f^{\ast}(u) = y \wedge (fi_{x})^{\ast}(u)$ for all $u \in \mathcal{O}(X)$.

Suppose $\bigvee\{y_{i} \mid i \in I\} = 1_{\mathcal{O}(Y)}$ be a cover of $Y$ such that each composite $fi_{y_{j}}: \downarrow y_{j} \rightarrow Y \rightarrow X$ is an open inclusion. For each $y \in  \mathcal{O}(Y)$, by lemma 4.2, we have $y \wedge y_{j} \lhd_{f} y \wedge y_{j}$ for all $j \in I$. Hence $y = \bigvee \{x \in \mathcal{O}(Y) \mid x  \lhd_{f} y\}$.  $\hfill\Box$

In general, a continuous posheaf is not a completely distributive posheaf. Indeed if we consider the case when $X$ be the terminal locale then a continuous posheaf is just a continuous lattice and  a completely distributive posheaf is just a completely distributive lattice.

\begin{definition} Let $f: Y \rightarrow X$ be a locale morphism. $y \in \mathcal{O}(Y)$ is said to be a $f$-coprime element if  $y \wedge f^{\ast}(u) \leq (x \vee z) \wedge f^{\ast}(u)$ implies that $y \wedge f^{\ast}(u) \leq x \wedge f^{\ast}(u)$ or $y \wedge f^{\ast}(u) \leq z \wedge f^{\ast}(u)$ for any $x, z \in \mathcal{O}(Y)$ and $u \in \mathcal{O}(X)$. Write $CP_{f}$ for the set of all $f$-coprime elements of $\mathcal{O}(Y)$. \end{definition}

 If $X$ being the terminal locale, then an $f$-coprime element of $\mathcal{O}(Y)$ is just a coprime element in $\mathcal{O}(Y)$. If $\mathcal{O}(Y)$ is the powerset lattice of a set then every singleton set is an $f$-coprime element. If $\mathcal{O}(Y) = [0, 1]$ be the unit interval then every element of $\mathcal{O}(Y)$ is  an $f$-coprime element.

\begin{proposition} Let $f: Y \rightarrow X$ be a locale morphism such that the corresponding frame sheaf of $f$ is a continuous posheaf. If $\mathcal{O}(Y)$ has enough $f$-coprime elements, i.e. $x = \bigvee \{y \leq x \mid y \in CP_{f}\}$ for every $x \in \mathcal{O}(Y)$ then the corresponding frame sheaf of $f$ is a completely distributive posheaf. \end{proposition}

{\bf Proof}\, We only need to show that $z \in CP_{f}$ and $z \ll_{f} x$ implies that $z \lhd_{f} x$. It is clear by the definition.  $\hfill\Box$

Recall that an internal frame $F$ in a localic topos $Sh(X)$ is called spatial if $F$ is isomorphic to a subframe of a power object. Classically we know every continuous frame,  moreover every completely distributive lattice is spatial under the assumption of AC. But we don't know weather it is still true in a localic topos. Under some additional conditions, we can show a partial result.

Let $f: Y \rightarrow X$ be a locale morphism such that the corresponding frame sheaf of $f$ is a completely distributive posheaf. We call $f: Y \rightarrow X$ an algebraic completely distributive posheaf if $y = \bigvee \{x \in \mathcal{O}(Y) \mid x  \lhd_{f} x \leq y\}$ for all $y \in \mathcal{O}(Y)$.

\begin{proposition} If $f: Y \rightarrow X$ is an algebraic completely distributive posheaf then $f: Y \rightarrow X$ is spatial. \end{proposition}

{\bf Proof}\, Write $\downdownarrows_{f} y = \bigvee \{x \in \mathcal{O}(Y) \mid x  \lhd_{f} x \leq y\}$ for  $y \in \mathcal{O}(Y)$. By Lemma 3.1, $G(u) = \downdownarrows_{f} f^{\ast}(u)$ defines a subsheaf of the frame sheaf $f: Y \rightarrow X$. We now show that the frame sheaf $f: Y \rightarrow X$ can be embed as a subsheaf into the powersheaf $\mathbb{P}G$.

Indeed, consider the morphism $\gamma: F \rightarrow \mathbb{P}G$ such that for each $u \in \mathcal{O}(X)$ and $y \leq f^{\ast}(u)$,
\[\gamma_{u}(y)(v) = \left\{\begin{array}{cc}
\downdownarrows_{f} (y\wedge f^{\ast}(v)), & v \leq f^{\sharp}(y)\\
\emptyset, &  v \not\leq f^{\sharp}(y)
\end{array}
\right.\]
Then it is straightforward to show that $\gamma$ is a frame monomorphism.

\end{document}


\newcommand{\Proof}{\noindent {\bf Proof.}}
\newcommand{\Section}[1]{
   \stepcounter{section}
   \bigskip\noindent
   {\bf\hbox{\thesection.~}#1}\par
   \nopagebreak
   \medskip
   \renewcommand{\theequation}{\thesection.\arabic{equation}}
   \setcounter{equation}{0}
   \setcounter{subsection}{0}}
\theoremstyle{definition}
\newtheorem{definition}{Definition}[section]
\newtheorem{remark}{Remark}[section]
\theoremstyle{plain}
\newtheorem{lemma}{Lemma}[section]
\newtheorem{theorem}{Theorem}[section]
\newtheorem{example}{Example}[section]
\newtheorem{corollary}{Corollary}[section]
\newtheorem{proposition}{Proposition}[section]
\newtheorem{question}{Question}[section]
\title{{\bf Partially Ordered Sheaves on a Locale. I}
\thanks{\it Project supported by NSFC (11171156) and PCSIRT (IRTO0742).}}
\author{{\small Wei He}\thanks{{{\it E-mail address:} weihe@njnu.edu.cn}}\\
 \small Institute
of Mathematics, Nanjing Normal University, Nanjing, 210097, China}

\date{}
\maketitle

\begin{abstract}
 In this paper, we investigate the order algebraic structure in the category of sheaves on a given locale
 $X$. Since every localic topos has a generating set formed by its
subterminal objects, we define a ``point" of a partially ordered
sheaf to be a morphism from a subterminal sheaf to the  partially
ordered sheaf. Using the concept of  ``points",  we investigate
the completeness of posheaves systemically. Some internal
characterizations of complete partially ordered sheaves and  frame
sheaves are given.  We also give
an explicit description of the construction of associated sheaf
locales and show directly that the category $Sh(X)$ of sheaves on
a locale $X$ is equivalent to the slice category $LH/X$ of locales
and local homeomorphisms over $X$. Applying this equivalence, we
give characterizations of partially ordered sheaves and complete
partially ordered sheaves in terms of sheaf locale respectively.

 \noindent{\bf Keywords:}  sheaf; partially ordered sheaf; sheaf locale.

\noindent{\bf Mathematics Subject Classifications}(2000):  18F20.
\end{abstract}

\section{Introduction}
The theory of partially ordered sets in a topos has been studied extensively (see [1], [2], [3]). When we restrict  our attention to a localic topos, i.e.
 the topos $Sh(X)$ of sheaves on a given locale $X$,  it is interesting to investigate the properties of ordered sheaves. In this paper,
 we first introduce the concept of partially ordered sheaves on a given locale $X$ which are partially ordered objects in the localic
topos $Sh(X)$. For every localic topos $Sh(X)$, it is well known
that it is not well-pointed, i.e. the terminal sheaf $1$ can not
generate $Sh(X)$. But the localic topos $Sh(X)$ has a generating set
formed by its subterminal objects. This implies that in a localic
topos, there exists enough ``points"  of morphisms from
subterminal sheaves to a sheaf $F$ such as $\hat{1} \rightarrow F$,
where $\hat{1}$ is a subsheaf of the terminal sheaf $1$. For a
given sheaf $F$, these points in $F$ will act similar as what points
in a set. By
 using the concept of ``points",  we  define the concepts such as upper
 bound and
 $supremum$ for a subsheaf of a partially ordered sheaf. This makes
 us can investigate those more complicated concepts such as complete
 partially ordered sheaves and frame sheaves.

The paper is organized in to five sections. In section 2, we
introduce the concept of partially ordered sheaves, and
investigate the basic properties of partially ordered sheaves.  In
section 3,  we investigate the completeness of partially ordered
sheaves, an internal characterization of complete partially
ordered sheaves is given. In particular, we give a
characterization of complete Heyting sheaves on $X$ which are
properly internal frames in the localic topos $Sh(X)$.  At the
end, we give an explicit description of the associated sheaf
locale of a given sheaf and show directly that the category
$Sh(X)$ of sheaves on a locale $X$ is equivalent to the slice
category $LH/X$ of locales and local homeomorphisms over $X$. Thus
characterizations  of partially ordered sheaf locales,  and
complete partially ordered sheaf locales are presented
respectively. Throughout this paper, when we write $X$ for a
locale, we will write $\mathcal{O}(X)$ for the corresponding
frame. Readers may refer to [5] for notations and terminology not
explicitly given here.

\section{Partially Ordered Sheaves}
\begin{definition} Let $X$ be a locale and $F$ a sheaf on $X$. $F$ is called a partially ordered sheaf (shortly posheaf)
if and only if $F$ satisfies the following conditions:

(POS1) $F(u)$ is a partially ordered set for every $u \in \mathcal{O}(X)$;

(POS2) every restriction map $F(u) \rightarrow F(v)$ for $v \leq u$ is order-preserving;

(POS3)  given two compatible families $\{s_{i} \in F(u_{i}) \mid i \in
 I\}$ and $\{t_{i} \in F(u_{i}) \mid i \in
 I\}$ with $s_{i} \leq t_{i}$ for any $i \in I$. If $\{s_{i} \in F(u_{i}) \mid i \in
 I\}$ patch to an element $s \in
 F(\bigvee u_{i})$ and  $\{t_{i} \in F(u_{i}) \mid i \in
 I\}$ patch to another element $t \in
 F(\bigvee u_{i})$, then $s \leq t$.\end{definition}

 Examples of partially ordered sheaves are numerous, for example the
 continuous
 real-valued function sheaf $C$ which sends each open $u \in
 \mathcal{O}(X)$ the set $C(u) = \{f: u
 \rightarrow \mathbb{R}\ \mid f\, \mbox{is continuous}\}$ of continuous real-valued functions on $u$ with pointwise order is a posheaf on $X$.
 But a sheaf of partially ordered sets in general does not
 satisfy the condition (POS3), hence not always a posheaf.

 \begin{example} Consider the continuous real-valued function sheaf $C$ on
 $X$. If we keep the pointwise  order on $C(u)$ for $u \not=
 1_{X}$ and take discrete order on
 $C(1_{X})$ then it is a sheaf of partially ordered sets
 but not a posheaf. \end{example}

 \begin{lemma} $F$ is a posheaf iff $F$ is a sheaf of partially ordered sets and satisfies (POS3).\end{lemma}

If $F$ is a posheaf on $X$, then the sheaf $F^{op}$ defined by $F^{op}(u) = F(u)^{op}$, the opposite poset of $F(u)$, and each restriction
map $F^{op}(u) \rightarrow F^{op}(v)$ is same as  $F(u) \rightarrow F(v)$, is a poshaf.   If $G$ is a subsheaf of $F$, it is clear
that $G$ is a posheaf for the induced order. Recall that an internal
partially ordered object  in a topos $\mathcal{\varepsilon}$ is an
object $A$ of $\mathcal{\varepsilon}$ with a subobject
 $\leq_{A} \rightarrowtail A \times A$ satisfying the following conditions:

  (i) (reflexivity) the diagonal $\delta: A \rightarrow A \times A$ can be factored through $\leq_{A} \rightarrowtail A \times A$;

  (ii) (antisymmetry) the intersection $\leq_{A} \cap \leq_{A}^{op}$ of subobjects is contained in the diagonal, where $\leq_{A}^{op}$ is the image of the composite
  $\leq_{A} \stackrel{l}{\rightarrowtail} A \times A \stackrel{t}{\rightarrow} A \times A$ with $ A \times A \stackrel{t}{\rightarrow} A \times A$ the twist map interchange the factors
  of the product;

  (iii) (transitivity) the subobject $\langle p_{1} v, p_{2} u\rangle: C \rightarrowtail A \times A$ can be factored through $\leq_{A}  \stackrel{l}{\rightarrowtail} A \times A$, where $C$ is defined as the following pullback with projections $u$ and $v$:

  \[
\xymatrix{C \ar[r]^{u}\ar[d]_{v}& \leq_{A}
\ar[d]^{p_1}\\
      \leq_{A}  \ar[r]^{p_{2}}& A  }
\]

 If $F$ is a posheaf, then the sub-presheaf $\leq_{F}: \mathcal{O}(X)^{op} \rightarrow Set$ of the product sheaf $F \times F$ defined by
 $\leq_{F}(u) = \{(x, y) \in F(u) \times F(u)\mid x \leq y\}$ for $u \in \mathcal{O}(X)$ is a sheaf and satisfying conditions (i)-(iii). Conversely, if we have a subsheaf
 $\leq_{F} \rightarrowtail F \times F$ satisfying conditions (i)-(iii), then $F$ is a posheaf.  So we have the following result.

 \begin{lemma}  $F$ is a  posheaf on a locale $X$ if and only if $F$ is an internal partially ordered object in the localic topos
 $Sh(X)$.\end{lemma}

 We know that a localic topos is in general not well-pointed, i.e. the terminal object $1$ is not a generator.  But every localic topos can be generated
 by the subobjects of its terminal object $1$. This implies that in a topos $Sh(X)$ of sheaves, those ``points"  $\hat{1} \rightarrow F$ of a sheaf $F$ can
 act somewhat as points in the category of sets. We define a point of a sheaf $F$ to be a  morphism $p: \hat{1} \rightarrow F$ with $\hat{1}$ a subsheaf of
 the terminal sheaf $1$. A point of the form $1 \rightarrow F$  will be called  a global point of $F$. For a point $p: \hat{1} \rightarrow F$ of $F$,
 we write $dom(p)$ for the largest open $u \in \mathcal{O}(X)$ with $p(u) \not= \emptyset$, i.e. $dom(u) = \bigvee \{u \in \mathcal{O}(X)
 \mid p(u) \not= \emptyset\}$, and call it the domain of $p$. If we look at the image of a point
 $p: \hat{1} \rightarrow F$,  a point $p$ of $F$ can
 be equivalently regarded as an element of $F(dom(p))$.  The set of all points of a sheaf $F$ will be denoted by $F_{p}$.

 \begin{definition} Let $F$ be a posheaf.  We define a partial order on the set $F_{p}$ as following
 \[p_{1} \leq p_{2}  \Leftrightarrow dom(p_{1}) \leq dom(p_{2}),\, {\mbox and}\,\,  p_{1}(dom(p_{1})) \leq p_{2}(dom(p_{2}))|_{dom(p_1)}\] where  $p_{2}(dom(p_{2}))|_{dom(p_1)}$ be the restriction of $p_{2}(dom(p_{2}))$ on $dom(p_1)$. \end{definition}

This definition of partial order is equivalent to saying that for
two points $p_{1}: \hat{1} \rightarrow F$ and  $p_{2}: \check{1}
\rightarrow F$, $ p_{1} \leq p_{2} $ if and only if there is a
 morphism $h: \hat{1} \rightarrow \check{1}$ such that
$\langle p_{1}, p_{2}h \rangle: \hat{1} \rightarrow F \times F$ can
be factored through $ \leq \rightarrowtail F \times F$.

 Given a
morphism $\alpha: F \rightarrow G$ of sheaves, we have a natural map
$F_{p} \rightarrow G_{p},  p\mapsto \alpha p$. For morphisms
$\alpha: F \rightarrow G$ and $\beta: F \rightarrow G$  of
posheaves, we define $\alpha \leq \beta$ if and only if $\alpha p
\leq \beta p$ for all points $p \in F_{p}$. Then we have

\begin{lemma}  For morphisms
$\alpha: F \rightarrow G$ and $\beta: F \rightarrow G$  of
posheaves, the following are equivalent:

(1)\, $\alpha \leq \beta$;

(2)\, $\forall u \in \mathcal{O}(X), \forall x \in F(u),
\alpha_{u}(x) \leq \beta_{u}(x)$ ;

(3)\, the diagonal $\langle \alpha, \beta \rangle: F \rightarrow G
\times G$ can be factored through $\leq_{G} \rightarrowtail G \times
G$. \end{lemma}

\begin{definition} Let $F$ and $G$ be partially ordered sheaves on a locale $X$, and $\alpha: F \rightarrow G$. $\alpha$ is called order-preserving if
for any two points $p_1, p_2$ of $F$, $p_{1} \leq p_{2}$ implies
$\alpha p_{1} \leq \alpha p_{2}$.\end{definition}

If $f: A \rightarrow B$ and $g: B \rightarrow C$ are
order-preserving, it is clear that the composite $gf: A \rightarrow
C$ is order-preserving.

\begin{proposition}  Let $F$ and $G$ be partially ordered sheaves on a locale $X$, and $\alpha: F \rightarrow G$ a morphism. The following conditions are equivalent:

(1) $\alpha$ is order-preserving;

(2)  $\forall u \in \mathcal{O}(X),\, \alpha_{u}: F(u) \rightarrow
G(U)$ is order-preserving;

(3)  the composite $\leq_{F} \rightarrowtail  F\times F
\stackrel{\alpha\times \alpha}{\rightarrow} G \times G$ can be factored
through $\leq_{G} \rightarrowtail G \times G$, i.e. we have a
morphism $\gamma: \leq_{F} \rightarrow \leq_{G}$ such that the
following square commutes

\[
\xymatrix{ F \times F\ar[r]^{\alpha\times \alpha}& G \times G\\
           \leq_{F}\ar[r]^{\gamma}\ar[u]^{}
           &\leq_{G}\ar[u]^{}}
\]
\end{proposition}

\begin{definition} Two partially ordered sheaves $F$ and $G$  on a locale
$X$ is said to be order isomorphic if there is a isomorphism
$\alpha: F \rightarrow G$ such that both $\alpha$ and its inverse
$\alpha^{-1}$ is order-preserving.\end{definition}

  Let $F$ be a posheaf on $X$
and $G$ a subsheaf of $F$. We call $G$ a downsheaf of $F$ if for any
two points $p,  p' \in F_{p}$,  $p \leq p'$ and $ p' \in G_{p}$
implies $p \in G_{p}$.

\begin{proposition} Let $F$ be a posheaf on $X$ and $G$ a subsheaf of
$F$. The following conditions are equivalent:

(1) $G$ is a downsheaf of $F$;

(2)  $\forall u \in \mathcal{O}(X),  G(u)$ is a downset of $F(u)$;

(3)  the classification map $\phi: F^{op} \rightarrow \Omega$ for
$G^{op}$ is order-preserving, where $\Omega$ be the
subobject classifier in $Sh(X)$.
\end{proposition}

{\bf Proof}\, The equivalence of (1) and (2) is clear.

(2) $\Rightarrow$ (3)\, Condition (2) is equivalently to say that $G^{op}(u)$ is an upper set of $F^{op}(u)$ for every $u \in \mathcal{O}(X)$.  For $u \in \mathcal{O}(X)$, $\phi_{u}: F^{op}(u) \rightarrow \Omega(u)$ is defined as
 $\phi_{u}(x) = \bigvee\{v \mid v \leq u, x|_{v} \in G(v)\}, x \in F(u)$.

 \[
\xymatrix{G(u) \ar[r]^{}\ar[d]_{}& 1(u)
\ar[d]^{}\\
      F^{op}(u)  \ar[r]^{\phi_{u}}& \Omega(u)}
\]
So if $x \leq^{op} y$ in $F^{op}(u)$, i.e. $y \leq x$ in $F(u)$, then $x|_{v} \in G(v)$ implies $y|_{v} \in G(V)$. Hence $\phi_{u}(x) \leq \phi_{u}(y)$.

(3) $\Rightarrow$ (2)\, If $x \leq y$ in $F(u)$ and $y \in G(u)$. Then $\phi_{u}(y) \leq \phi_{u}(x)$ and $\phi_{u}(y)= u$. Hence $\phi_{u}(x) = u$ this implies $x \in G(u)$. $\hfill\Box$

Let $F$ be a posheaf on $X$ and $p$ a point of $F$. We have a
downsheaf  $\downarrow p$ generated by $p$:

\[\downarrow p(u) = \left\{\begin{array}{cc}
\{x \in F(u)\mid x \leq p(u)\}, & p(u) \not= \emptyset\\
\emptyset, & p(u) =\emptyset
\end{array}
\right.\]

The downsheaf of the form $\downarrow p$ will be called a principle
ideal of the posheaf $F$.

 Dually, we define an uppersheaf of $F$ to be a downsheaf of $F^{op}$. A special class of uppersheaves for which we will call principe filters has the form
 \[\uparrow p(u) = \left\{\begin{array}{cc}
\{x \in F(u)\mid  p(u) \leq x\}, & p(u) \not= \emptyset\\
\emptyset, & p(u) =\emptyset
\end{array}
\right.\] where $p$ is a point of $F$.

Let  $F$  be a sheaf on a locale $X$ and $u \in \mathcal{O}(X)$. We
have a sheaf $F^{u}$ on the open $u$ such that $F^{u}(v) = F(v)$ for
every $v \leq u$.
  We call $F^{u}$ the restriction of $F$ to $u$. We some time also regard $F^{u}$ as a subsheaf of $F$ for which we take $F^{u}(w) = \emptyset$
  for any $w \not \leq u$.

   Recall the powersheaf $\mathbb{P}F$ of $F$ defined by  $\mathbb{P}F(u) = Sub(F^{u})$ with restriction maps
   $\mathbb{P}F(u) \rightarrow \mathbb{P}F(v), S \mapsto S^{v}$ for every $v \leq u$ in $\mathcal{O}(X)$.  $\mathbb{P}F$ is a partially ordered sheaf for
   the subsheaves inclusion order. For every element
    $x \in F(u)$, we regard $x$ as a point of $F$ and thus have a downsheaf $\downarrow x$. Hence we have a morphism
    \[\downarrow: F \rightarrow \mathbb{P}F\] such that for every $u \in \mathcal{O}(X)$ and $x \in F(u)$, $\downarrow(u)(x) = \downarrow
    x$. We call $\downarrow: F \rightarrow \mathbb{P}F$ the principle ideal
    embedding. Similarly, we can construct a down-powersheaf
    $\mathbb{D}F$ of $F$ such that  $\mathbb{D}F(u) =
    Dow(F^{u})$ where $Dow(F^{u})$ is the set of all downsheaves of
    $F^{u}$ with the same restriction maps with $\mathbb{P}F$. Then
    $\mathbb{D}F$ is a subsheaf of $\mathbb{P}F$, so
we have an inclusion $\mathbb{D}F \rightarrowtail \mathbb{P}F$ and
the principle ideal embedding $\downarrow: F \rightarrow
\mathbb{P}F$ can be factored through $\mathbb{D}F$.

Now we consider the generalization of the important concept of
Galois connection in classical order theory.

\begin{definition} For order-preserving $\alpha: F
\rightarrow G$ and $\beta: G \rightarrow F$, $\alpha$ is left
adjoint to $\beta$ ($\beta$ is right adjoint to $\alpha$), written
$\alpha \dashv \beta$ if and only if the relations $\alpha x \leq y$ and $x \leq \beta y$ are equivalent for all points $x \in F_{p}$ and
$y \in G_{p}$.
\end{definition}

As in the standard case, an adjoint, if it exists, is uniquely defined and the adjointness can be characterized by  $\alpha \dashv \beta$
if and only if $1_{F} \leq \beta\alpha$ and $\alpha\beta \leq
1_{G}$.

Similar to the case in the category of sets, an order-preserving morphism may has neither a left adjoint nor a
right adjoint. But if an order-preserving morphism has right (left)
adjoint then it is unique.

\begin{lemma} $(\alpha \dashv \beta):  F
\rightarrow G$ if and only if $\alpha_{u}: F(u) \rightarrow G(u)$ is left adjoint to $\beta_{u}: G(u) \rightarrow F(u)$ for every
$u \in \mathcal{O}(X)$.\end{lemma}

Now we construct a left adjoint for the inclusion  $\mathbb{D}F
\rightarrowtail \mathbb{P}F$.

If $S \in Sub(F)$ is a subsheaf of $F$, we write $\downarrow S$
the sub-preshaef of $F$ such that $\downarrow S(u) = \{x \in F(u)
\mid \exists \{u_{i} \mid i\in I\} \subset \mathcal{O}(X), x_{i} \in
S(u_{i}),  u = \bigvee u_{i}, x|_{u_{i}} \leq x_{i}, i\in I\}$ for
$u \in \mathcal{O}(X)$. Then $\downarrow S$ is a downsheaf of $F$.
Thus we have a morphism $\downarrow ( ): \mathbb{P}F \rightarrow
\mathbb{D}F$ such that for every $u \in \mathcal{O}(X)$ and $S \in
\mathbb{P}F(u)$, $\downarrow (u)(S) = \downarrow
    S$.

\begin{proposition} $\downarrow( ): \mathbb{P}F \rightarrow \mathbb{D}F$ is
order-preserving which is left adjoint to the inclusion $\mathbb{D}F
\rightarrowtail \mathbb{P}F$, i.e. $\downarrow S \subset G
\Leftrightarrow S \subset G$ for any $S \in Sub(F)$ and $G \in
Dow(F)$.
\end{proposition}

\begin{corollary} Let $F$ be a posheaf and $G$ a subsheaf of $F$.
$G$ is a downsheaf iff $G = \downarrow G$.\end{corollary}

\section{Complete Posheaves}

Let $F$ be a posheaf on a locale $X$ and $A$ a subsheaf of $F$. A
 point $x \in F_{p}$ is called an upper bound of $A$ if for any point $y
 \in A_{p}$ we have $y \leq x$. A least upper bound, otherwise known as a $supremum$ or even $sup$  for $A$,  is
 an upper bound $p$ for $A$ such that any upper bound $r$ for $A$ we
 have $p \leq r$. The $supremum$ for $A$, if it exists, is unique. We
 will write $\bigvee A$ for the $supremum$ of $A$. Dually, we define a lower bound $l$ of $A$ to be an upper bound of $A$
 in the opposite posheaf $F^{op}$,  and $infimum$ (shortly $inf$) for $A$ to be the largest lower bound of $A$, i.e. the least upper bound in $F^{op}$ and
 write it $\bigwedge A$.

 Suppose $B$ is a sub-presheaf of a posheaf $F$, and $\hat{B}$ is the subsheaf of $F$ generated by $B$. Then it is clear that $\hat{B}$ and
 $B$ have the same upper bound set and the same lower bound set. If $\alpha: F \rightarrow G$ is an order-preserving morphism and $S$ is a subsheaf of $F$,
 we call the sub-presheaf $\alpha S$ of $G$ defined by $\alpha S(u) = \alpha_{u}(S(u))$ the preimage of $S$ under $\alpha$,
 and the subsheaf $\hat{\alpha S}$ generated by $\alpha S$ the image of $S$ under $\alpha$.  Suppose $\bigvee S$ exists for
 a subsheaf $S$ of $F$, we say that $\alpha$ preserves $\bigvee S$ if $\bigvee \hat{\alpha S}$ exists and  $\alpha \bigvee S = \bigvee \hat{\alpha S}$ holds.
  Similarly we define $\alpha$ preserves $\bigwedge S$ if $\bigwedge \hat{\alpha S}$ exists and  $\alpha \bigwedge S = \bigwedge \hat{\alpha S}$ holds.

 \begin{proposition} If $(\alpha \dashv \beta):  F
\rightarrow G$ then $\alpha$ preserves any suprema that exist in $F$. Dually, $\beta$ preserves any infima that exist in $G$. \end{proposition}

{\bf Proof}\, Let $S \subseteq F$ be a subsheaf of $F$ with $\bigvee
S$ exists. It is clear that $\alpha \bigvee S$ is an upper bound of
$\alpha S$ since $\alpha$ is order preserving. Suppose $y \in G_{p}$ is another upper bound of $\alpha
S$. Then for any point $x$ in $S$, $\alpha x \leq y $,  hence $ x \leq \beta y$ by
adjointness. Thus we have $\bigvee S \leq
\beta y $  which means  $ \alpha \bigvee S \leq y$.  This shows that  $\bigvee
\hat{\alpha S} = \alpha \bigvee S$. The left follows from the fact
that $(\alpha \dashv \beta):  F \rightarrow G$ implies $(\beta
\dashv \alpha): G^{op} \rightarrow F^{op}$. $\hfill\Box$

\begin{definition} Let $F$  be a posheaf on a locale $X$. $F$ is called complete if the principle ideals embedding
$\downarrow: F \rightarrow \mathbb{D}F$  has a left
adjoint.\end{definition}

 \begin{lemma} Let $F$ be a posheaf on a locale $X$. If each restriction map $F(u) \rightarrow F(v)$ is surjective and it has a left adjoint, then for
each $u, v \in \mathcal{O}(X)$, the following square commutes:
  \[
\xymatrix{ F(u \vee v)\ar[r]^{}&  F(u)\\
           F(v)\ar[r]^{}\ar[u]^{l_{v u\vee v}}
           & F(u \wedge v)\ar[u]^{l_{u\wedge v u}}}
\] where $l_{v u\vee v}: F(v) \rightarrow F(u \vee v)$ and $l_{u\wedge v u}: F(u \wedge v) \rightarrow F(u)$ are the left adjoint of the restriction maps $F(u \vee v) \rightarrow
F(v)$ and $F(u) \rightarrow F(u \wedge v)$ respectively.\end{lemma}

{\bf Proof}\, We note that the above square always commutes for $u \leq v$ or $v \leq u$.
Suppose $v \not\leq u$, $u \not\leq v$ and let $x \in F(v)$.  We first have $((l_{v u\vee v}(x))|_{u})|_{u\wedge v} = ((l_{v u\vee v}(x))|_{v})_{u\wedge v} = x|_{u\wedge v}$. Suppose $y \in F(u)$ and $y|_{u\wedge v} = x|_{u\wedge v}$,
we have an unique element $z \in F(u \vee v)$ such that $z|_{v} = x,  z|_{u} = y$. Thus $l_{v u\vee v}(x) \leq z$ and $(l_{v u\vee v}(x))|_{u} \leq z|_{u} = y$. This shows that $(l_{v u\vee v}(x))|_{u}$ is the least element $t$ in $F(u)$ such that $t|_{u\wedge v} = x|_{u\wedge v}$. Hence $l_{u\wedge v u}(x|_{u\wedge v}) = (l_{v u\vee v}(x))|_{u}$.
 $\hfill\Box$

\begin{proposition} Let $F$  be a posheaf on a locale $X$. The following are equivalent:

(1) $F$ is complete.

(2) For every downsheaf $S$ of $F$, $\bigvee S$  exists and it can
be extended to a global point $p$ of $F$ such that for each $u \in
\mathcal{O}(X)$, $p|_{u}$ is the least  element in $F(u)$ satisfying
$S^{u} \subseteq \downarrow p|_{u}$.

(3)  For every subsheaf $S$ of $F$, $\bigvee S$  exists and it can
be extended to a global point $p$ of $F$ such that for each $u \in
\mathcal{O}(X)$, $p|_{u}$ is the least element in $F(u)$ satisfying
$S^{u} \subseteq \downarrow p|_{u}$.

(4) Every $F(u)$ is a complete lattice,  and each restriction map $F(u) \rightarrow F(v)$ is surjective and it has both a left adjoint and
a right adjoint.
\end{proposition}

{\bf Proof}\,  $(1) \Leftrightarrow (2)$ The statement that the
principle ideals embedding $\downarrow: F \rightarrow \mathbb{D}F$
has a left adjoint is equivalent to say that we have a morphism
$sup: \mathbb{D}F \rightarrow F$ such that for any downsheaf $S$ of
$F$, $sup_{u}(S^{u})$ is the least element in $F(u)$ satisfying
$S^{u} \subseteq \downarrow sup_{u}(S^{u})$ for $u \in \mathcal{O}(X)$.
This is equivalent to the condition that for any downsheaf $S$ of
$F$, $\bigvee S$ exists and it can be extended to a global point
$sup_{1_{X}}(S) \in F(1_{X})$ such that for any $v \in
\mathcal{O}(X)$, $sup_{1_{X}}(S)|_{v}$ is the least upper bound of
$S^{v}$, i.e. the following square commutes
 \[
\xymatrix{Dow(F) \ar[r]^{sup_{1_{X}}}\ar[d]_{}& F(1_{X})
\ar[d]^{}\\
      Dow(F^{v})  \ar[r]^{sup_{v}}& F(v)}
\]

$(3) \Rightarrow (2)$ is clear.  To show $(2) \Rightarrow (3)$, it suffice to show that for any subsheaf $S$
of $F$, $\bigvee \downarrow S$ exists implies that $\bigvee S$ exists  and $\bigvee S = \bigvee  \downarrow S$. But it is clear since $ \downarrow S$ and $S$ has same
upper bound set.

 $(1) \Rightarrow (4)$    Suppose $x \in F(u)$, we regard $x$ as a point of $F$
 and thus a subsheaf of $F$. Applied $x$ to the above commute square we have an element $sup_{1_{X}}(x) \in F(1_{X})$ such that $sup_{1_{X}}(x)|_{u} = x$. Hence
 the restriction map $F(1_{X}) \rightarrow F(u)$  is surjective.
If $F$ is a complete posheaf, then it has a largest point- the
$supremum$ of $F$. This implies that each $F(u)$ has a largest
element,
 and the restriction maps preserve
 the top elements. $F$ also has a least point- the  $supremum$ of the leat subsheaf of $F$. This means that each $F(u)$ has a bottom element and
 the restriction maps preserve
 the bottom elements. Suppose $A \subseteq F(u)$ for an $u \in \mathcal{O}(X)$, consider the presheaf $\bar{A}$ determined by
  $A$: \[\bar{A}(u) = \left\{\begin{array}{cc}
\{x \in F(v)\mid \exists y \in A,  x = y|_{v}\}, & v \leq u\\
\emptyset, & v \not\leq u
\end{array}
\right.\]  Write $\hat{A}$ the subsheaf generated by $\bar{A}$, then $\bigvee \hat{A}$ exists. It is clear that $\bigvee \hat{A}$ is also the least upper
 bound of $A$ in the poset $F(u)$. Hence $F(u)$ is complete. For $v \leq u$ in $\mathcal{O}(X)$,  the commutative  of the square
 \[
\xymatrix{\mathbb{D}F(u) \ar[r]^{sup_{u}}\ar[d]_{}& F(u)
\ar[d]^{}\\
      \mathbb{D}F(v)  \ar[r]^{sup_{v}}& F(v)}
\] implies that $\bigvee \hat{A}|_{v} = \bigvee \hat{A}^{v}$,
i.e. $(\bigvee A)|_{v} = \bigvee (A|_{v})$. Hence the restriction map
$F(u) \rightarrow F(v)$ preserves joins, so has a right adjoint. To
show the restriction map $F(u) \rightarrow F(v)$ has a left adjoint,
we note that for any element $x \in F(v)$, if we regard $x$ as a
point and thus a subsheaf of $F$, there exists a least element $y
\in F(u)$ such that $x \leq y|_{v}$, i.e. $min \{y \in F(u) \mid x
\leq y|_{v}\}$ exists. Hence the restriction map $F(u) \rightarrow
F(v)$ has a left adjoint.

 $(4) \Rightarrow (3)$  Suppose $S$ is a subsheaf of $F$, write $u = \bigvee \{v \in \mathcal{O}(X) \mid S(v) \not= 0\}$.
 For $v \leq u$ in $\mathcal{O}(X)$, denote $l_{vu}: F(v) \rightarrow F(u)$ the left adjoint of the restriction map $F(u) \rightarrow F(v)$. Let
  $s = \bigvee \{l_{vu}(\bigvee S(v)) \mid v \leq u\}$, then it is clear that $s$ is the least upper bound of $S$. Write
  $\bar{s} = l_{u1_{X}}(s)$, then $\bar{s}|_{u} = s$.  We first show $\bar{s}|_{v} = s|_{v} = \bigvee
  S^{v}$ for  $v \leq u$.

  In fact $s|_{v} = (\bigvee \{(l_{wu}(\bigvee S(w)))|_{v} \mid w \leq v\}) \vee (\bigvee \{(l_{wu}(\bigvee S(w)))|_{v} \mid w \leq u,  w \not\leq v\})
  = (\bigvee \{(l_{wv}(\bigvee S(w)))|_{v} \mid w \leq v\}) \vee (\bigvee \{(l_{wu}(\bigvee S(w)))|_{v} \mid w \leq u,  w \not\leq v\})
  = (\bigvee S^{v}) \vee (\bigvee \{(l_{wu}(\bigvee S(w)))|_{v} \mid w \leq u,  w \not\leq v\})$. For each $w \leq u,  w \not\leq v$, by lemma 3.1,
   we have $(l_{wu}(\bigvee S(w)))|_{v} = (l_{w v\vee w}(\bigvee S(w)))|_{v} = l_{v\wedge w v}((\bigvee S(w))|_{v \wedge w}) \leq
   l_{v\wedge w v}(\bigvee S(v \wedge w))$ since every restriction map preserves sups. Hence
   $\bigvee \{(l_{wu}(\bigvee S(w)))|_{v} \mid w \leq u,  w \not\leq v\} \leq \bigvee S^{v}$. This shows $ s|_{v} = \bigvee
  S^{v}$.

  Suppose $v \not\leq u$. By lemma 3.1 and the above result,  we have  $\bigvee S^{v} = l_{u\wedge v v}(\bigvee S^{u \wedge v}) = l_{u\wedge v v}(s|_{u \wedge v}) = (l_{u u\vee v}(s))|_{v}
  = \bar{s}|_{v}$.
 $\hfill\Box$

\begin{example} Let $F$ be a sheaf on locale $X$. Consider the  power
sheaf $\mathbb{P}F$ of $F$. For any $u \in \mathcal{O} (X)$,
$\mathbb{P}F(u) = sub(F^{u})$ is a complete lattice.  Moreover, the left adjoint $l_{vu}: F(v) \rightarrow F(u)$
of the restriction map $F(u) \rightarrow F(v)$ sends each subsheaf $S$ of $F^{v}$ to a subsheaf $\check{S}$ of $F^{u}$ defined by for any $w \leq u$,
\[\check{S}(w) = \left\{\begin{array}{cc}
S(w), & w \leq v\\
\emptyset, & w \not\leq v
\end{array}
\right.\]
 The right adjoint of  the restriction map $F(u) \rightarrow F(v)$ sends each subsheaf $S$ of $F^{v}$ to a subsheaf $\hat{S}$ of $F^{u}$
 generated by the presheaf $\bar{S}$ with $\bar{S}(w) = \{x \in F(w)\mid x|_{w\wedge v} \in S(w \wedge v)\}$, $w \leq u$.
 Hence
$\mathbb{P}F$ is a complete partially ordered sheaf. Similarly, we can show the down-powersheaf
    $\mathbb{D}F$ of $F$ is a complete partially ordered sheaf. \end{example}

\begin{example} Let $X$ be a topological space. Consider the sheaf
$LSC_{X}$ of lower semicontinuous functions into the unit interval
$[0,1]$ on $X$. Then each $LSC_{X}(u)$ is a complete lattice, and
each restriction map $LSC_{X}(u) \rightarrow LSC_{X}(v)$ preserves
joins since joins are pointwise. For opens $v \leq u$, and a lower
semicontinuous map $f: v \rightarrow [0, 1]$, define $\bar{f}: u
\rightarrow [0, 1]$
 \[\bar{f}(x) = \left\{\begin{array}{cc}
f(x), & x \in v\\
0, & x \in u\setminus v
\end{array}
\right.\] Then $\bar{f}$ is lower semicontinuous and the
corresponding $f \mapsto \bar{f}$ forms a left adjoint of the
restriction map  $LSC_{X}(u) \rightarrow LSC_{X}(v)$.  Hence
 $LSC_{X}$ is a complete partially ordered sheaf.
\end{example}

\begin{corollary}  Let $F$  be a posheaf on a locale $X$.  $F$ is complete if and only if $F^{op}$ is complete. \end{corollary}

We write $SCL$ for the category of all complete lattices and
surjective maps preserving arbitrary sups and arbitrary infs.

\begin{corollary} Let $X$ be a locale. $F$ is a complete posheaf on $X$ if and only if $F$ is a sheaf over $SCL$  satisfying the condition
(POS3).  \end{corollary}

Let $F$ and $G$ be complete posheaves on a locale $X$ and let
$\alpha: F \rightarrow G$ be an order-preserving morphism. Then we
have an order preserving morphism $\alpha_{\ast}: \mathbb{P}F
\rightarrow \mathbb{P}G$ such that for each $u \in \mathcal{O}(X)$
and  $S \in sub(F^{u})$, ${\alpha_{\ast}}_{u}: sub(F^{u})
\rightarrow sub(G^{u})$ maps $S$ to the image of $S$ under $\alpha$.
We call $\alpha$ an sup-preserving morphism if the following square
commutes:
 \[
\xymatrix{\mathbb{P}F \ar[r]^{\alpha_{\ast}}\ar[d]_{sup_{F}}&
\mathbb{P}G
\ar[d]^{sup_{G}}\\
      F  \ar[r]^{\alpha}& G}
\] where $sup_{F}: \mathbb{P}(F) \rightarrow F$ and $sup_{G}: \mathbb{P}(G) \rightarrow
G$ are  the left adjoint of the principle ideal embedding
$\mathbb{P}F \rightarrow F$ and $\mathbb{P}G \rightarrow G$
respectively.

\begin{proposition}  Let $F$ and $G$ be complete posheaves on a locale $X$ and $\alpha: F
\rightarrow G$ be an order-preserving morphism. The following
conditions are equivalent:

(1)  $\alpha: F \rightarrow G$  is an sup-preserving morphism.

(2) For each $u \in \mathcal{O}(X)$, $\alpha_{u}: F(u) \rightarrow
G(u)$ preserves joins and the following square commutes for any $v
\leq u$ in $\mathcal{O}(X)$
 \[
\xymatrix{F(u)\ar[r]^{\alpha_{u}}&  G(u)\\
           F(v)\ar[r]^{\alpha_{v}}\ar[u]^{f_{uv}}
           & G(v)\ar[u]^{g_{uv}}}
\] where $f_{uv}: F(v) \rightarrow F(u)$ and $g_{uv}: G(v) \rightarrow G(u)$ are the left adjoint of the restriction maps $F(u) \rightarrow
F(v)$ and $G(u) \rightarrow G(v)$ respectively.

(3) $\alpha$ has a right adjoint, i.e. there exists an adjoint pair $\alpha\vdash\beta: F \rightarrow G$.
\end{proposition}

{\bf Proof}\,  $(1) \Rightarrow (2)$ Suppose $u \in \mathcal{O}(X)$
and  $A \subseteq F(u)$. Write $\hat{A}$ the subsheaf of $F$
generated by $A$, $\hat{\alpha A}$ the image of $\hat{A}$ under
$\alpha$. Then $\bigvee \hat{A} = \bigvee A$ and $\bigvee
\hat{\alpha A} = \bigvee \alpha A$. Hence we have
$\alpha_{u}(\bigvee A) = \bigvee \alpha_{u}(A)$ by applying the
commutative square of the definition. For each $v \leq u$ in
$\mathcal{O}(X)$ and  $x \in F(v)$, we regard $x$ as a point hence a
shbsheaf of $F^{u}$. We have $sup_{F u}(x) = f_{uv}(x),  sup_{G
u}(x) = g_{uv}(x)$. Hence $\alpha_{u} f_{uv}(x) = g_{uv}
\alpha_{v}(x)$ by the definition.

$(2) \Rightarrow (1)$ Suppose $u \in \mathcal{O}(X)$,  $S$ is a
subsheaf of $F^u$ and $\hat{\alpha S}$ the image of $S$ under
$\alpha$. Write $w = \bigvee \{v \in \mathcal{O}(X) \mid S(v) \not=
\emptyset \}$. Then $\bigvee S = \bigvee \{f_{wv}(\bigvee S(v)) \mid
v \leq w \}$. Thus $\bigvee \hat{\alpha S} = \bigvee
\{g_{wv}(\bigvee \alpha_{v} (S(v))) \mid v \leq w \} = \bigvee
\{g_{wv} \alpha_{v} (\bigvee S(v)) \mid v \leq w \} = \bigvee
\{\alpha_{w} f_{wv}(\bigvee S(v)) \mid v \leq w \} = \alpha_{v}
(\bigvee \{f_{wv}(\bigvee S(v)) \mid v \leq w \} = \alpha \bigvee
S$. By the completeness of $F$, we know $sup_{Fu}(S) = f_{uw}
(\bigvee S)$ and $sup_{Gu}(\hat{\alpha S}) = g_{uw}(\bigvee
\hat{\alpha S})$. Hence by the commutative square in (2), we have
$\alpha_{u} sup_{Fu}(S) = sup_{Gu}(\hat{\alpha S})$. This shows that
$\alpha$ is an sup-preserving morphism.

$(2) \Leftrightarrow (3)$ By lemma 2.6, $\alpha$ has a right adjoint $\beta$ if and only if  $\alpha_{u}: F(u) \rightarrow G(u)$ has a right adjoint $\beta_{u}: G(u) \rightarrow F(u)$ for every $u \in  \mathcal{O}(X)$ and the following square commutes for any $v
\leq u$ in $\mathcal{O}(X)$
\[
\xymatrix{G(u) \ar[r]^{\beta_{u}}\ar[d]_{}&
F(u)
\ar[d]^{}\\
      G(v)  \ar[r]^{\beta_{v}}& F(v)}
\] It equivalents to that each  $\alpha_{u}: F(u) \rightarrow G(u)$ preserves joins and  the following square commutes for any $v
\leq u$ in $\mathcal{O}(X)$
 \[
\xymatrix{F(u)\ar[r]^{\alpha_{u}}&  G(u)\\
           F(v)\ar[r]^{\alpha_{v}}\ar[u]^{f_{uv}}
           & G(v)\ar[u]^{g_{uv}}}
\] by the uniqueness of adjoint.
$\hfill\Box$

\begin{example} Let $F$ and $G$ be sheaves on a locale $X$, and
$\alpha: F \rightarrow G$ be a morphism. Consider the image morphism
$\alpha_{\ast}: \mathbb{P}F \rightarrow \mathbb{P}G$. For each $u
\in \mathcal{O}(X)$,  $\alpha_{\ast}(u): \mathbb{P}F(u) \rightarrow
\mathbb{P}G(u)$ preserves  joins of subsheaves.   The left adjoint $
\mathbb{P}F(v) \rightarrow \mathbb{P}F(u)$ for a restriction map
$\mathbb{P}F(u) \rightarrow \mathbb{P}F(v)$ sends each $S \in F^{v}$
to its minimal extension $\tilde{S} \in F^{u}$ defined by
\[\tilde{S}(w) = \left\{\begin{array}{cc}
S(w), & w \leq v\\
\emptyset, & w \not\leq v
\end{array}
\right.\] Hence the square in proposition 3.3 commutes. This shows
that $\alpha_{\ast}: \mathbb{P}F \rightarrow \mathbb{P}G$ is an
sup-preserving morphism. \end{example}

Now we consider the finite completeness of posheaves.

\begin{definition} Let $F$ be a posheaf on $X$.  $F$ is said to be finite sup-complete if  $F \rightarrow 1 $ and the diagonal $F \rightarrow F\times F$ both has a left adjoint.\end{definition}

\begin{proposition} Let $F$ be a posheaf on $X$. The followings are equivalent:

(1) $F$ is finite sup-complete.

(2) For every  $u \in \mathcal{O}(X)$, $F(u)$ is an sup-semilattice and every restriction map $F(u) \rightarrow F(v)$ preserves finite joins. \end{proposition}

Dually, we define a posheaf $F$ to be finite inf-complete if and only if  $F \rightarrow 1 $ and the diagonal $F \rightarrow F\times F$ both has a right adjoint.

\begin{proposition} Let $F$ be a posheaf on $X$. The followings are equivalent:

(1) $F$ is finite inf-complete.

(2) For every  $u \in \mathcal{O}(X)$, $F(u)$ is a inf-semilattice and every restriction map $F(u) \rightarrow F(v)$ preserves finite meets. \end{proposition}

\begin{definition} Let $F$ be a posheaf on $X$.  $F$ is said to be finite complete if it is both finite sup-complete and finite inf-complete \end{definition}

\begin{proposition} Let $F$ be a posheaf on $X$. The followings are equivalent:

(1) $F$ is finite complete.

(2) For every  $u \in \mathcal{O}(X)$, $F(u)$ is a lattice and every restriction map $F(u) \rightarrow F(v)$ preserves finite joins and finite meets. \end{proposition}

By proposition 3.2, we know that every complete posheaf is finite
complete.

Now we consider the generalization of another very important class
of order algebras- complete Heyting algebras (or frames).  For every
complete posheaf $F$, we define a meet morphism $\mu_{F}: F \times
\mathbb{P}F \rightarrow \mathbb{P}F$ of points with subsheaves such that for
each $u \in \mathcal{O}(X)$, $x \in F(u)$, and $S \in sub(F^{u})$,
$\mu_{F}(u)(x, S)$ be the subsheaf of $F^{u}$ generated by the
sub-presheaf $\bar{S}$ of $F^{u}$ with $\bar{S}(v) = \{x|_{v} \wedge
y \mid y \in S(v)\}$ for each $v \leq u$. Note that since each
restriction map preserves meets so $\bar{S}$ is indeed a
sub-presheaf of $F^{u}$.

\begin{definition} Let $F$ be a complete posheaf on a locale $X$.
$F$ is  said to be a complete Heyting sheaf (or frame sheaf) if the
following square commutes:

 \[
\xymatrix{F \times \mathbb{P}F \ar[r]^{1_{F}\times
sup_{F}}\ar[d]_{\mu_{F}}&  F \times F
\ar[d]^{m_{F}}\\
      \mathbb{P}F   \ar[r]^{sup_{F}}& F}
\]
 where $m_{F}: F \times F \rightarrow F$ is the
right adjoint of the diagonal $F \rightarrow F \times F$.
\end{definition}

\begin{proposition}  Let $F$ be a complete posheaf on a locale $X$.
The following conditions are equivalent:

(1) $F$ is a complete Heyting sheaf.

(2) Each $F(u)$ is a complete Heyting algebra for $u \in
\mathcal{O}(X)$, and if $v \leq u$ in $\mathcal{O}(X)$,  then $x
\wedge l_{uv}(y) = l_{uv}( x|_{v} \wedge y)$ holds for any $x \in
F(u), y \in F(v)$ where $l_{uv}: F(v) \rightarrow F(u)$ is the left
adjoint of the restriction map $F(u) \rightarrow F(v)$.
\end{proposition}

{\bf Proof}\, $(1) \Rightarrow (2)$  Suppose $x \in F(u), S
\subseteq F(u)$,  $u \in \mathcal{O}(X)$. Write $\hat{S}$ for the
subsheaf of $F^{u}$ generated by $S$. Then we have $\bigvee \hat{S}
=\bigvee S$, and $\bigvee \mu_{F}(u)(x, S) = \bigvee \{x \wedge s
\mid s \in S\}$. By the definition, we have $x \wedge \bigvee S =
\bigvee \{x \wedge s \mid s \in S\}$. For $v \leq u$ in
$\mathcal{O}(X)$, let $x \in F(u), y \in F(v)$. Regard $y$ as a
point and thus a subsheaf of $F^{u}$,  we have $sup_{F}(u)(y) =
l_{uv}(y)$ and $sup_{F}(u)(\mu_{F}(x, y)) = l_{uv}(x|_{v} \wedge
y)$. Hence $x \wedge l_{uv}(y) = l_{uv}( x|_{v} \wedge y)$ by the
definition.

$(2) \Rightarrow (1)$  Suppose  $x \in F(u)$ and $ S \in sub(F^u)$ a
subsheaf of $F^{u}$, $u \in \mathcal{O}(X)$. Write $w = \bigvee \{v
\in \mathcal{O}(X) \mid S(v) \not= \emptyset\}$. Then $\bigvee S =
\bigvee \{l_{wv}(\bigvee S(v)) \mid v \leq w\}$, and $sup_{F}(u)(S)
= l_{uw}(\bigvee S)$. Thus $sup_{F}(u)(\mu_{F}(u)(x, S)) =
l_{uw}(\bigvee \mu_{F}(u)(x, S)) = l_{uw}(\bigvee
\{l_{wv}(\bigvee\{x|_{v} \wedge y \mid y\in S(v)\} \mid v \leq w\})
= l_{uw}(\bigvee\{l_{wv}(x|_{v} \wedge \bigvee S(v)) \mid v \leq
w\}) = l_{uw}(\bigvee\{x|_{w} \wedge l_{wv}(\bigvee S(v)) \mid v
\leq w\} = l_{uw}(x|_{w} \wedge \bigvee \{l_{wv}(S(v)) \mid v \leq
w\}) = x \wedge l_{uw}(\bigvee \{l_{wv}(S(v)) \mid v \leq w\}) = x
\wedge sup_{F}(u)(S)$. $\hfill\Box$

\begin{example}  Every  power
sheaf $\mathbb{P}F$ of $F$ is a complete Heyting sheaf.  Moreover, every down-powersheaf
    $\mathbb{D}F$ of $F$  is a complete Heyting sheaf.  \end{example}

  Let $F$ and $G$ be two finite inf-complete posheaves on  $X$, we call a morphism $\alpha: F \rightarrow G$ of sheaves preserving finite meets if the following square commutes:

 \[
\xymatrix{F \times F \ar[r]^{\alpha\times
\alpha}\ar[d]_{m_{F}}&  G \times G
\ar[d]^{m_{G}}\\
      F   \ar[r]^{\alpha}& G}
\] where $m_{F}: F \times F \rightarrow F$ and $m_{G}: G \times G \rightarrow G$ represent the right adjoint of the diagonals  $F \rightarrow F \times F$ and
$G \rightarrow G \times G$ respectively.

\begin{lemma}  $\alpha: F \rightarrow G$ preserving finite meets if and only if for each $u \in \mathcal{O}(X)$, $\alpha_{u}: F(u) \rightarrow G(u)$
 preserves finite meets. \end{lemma}

\begin{definition} Let $F$ and $G$ be two frame sheaves on $X$.  A  morphism $\alpha: F \rightarrow G$ of sheaves is said to be a frame morphism if
$\alpha$ is a sup-preserving morphism which also preserves finite meets.  \end{definition}

\begin{lemma}  Let $F$ and $G$ be two frame sheaves on $X$ and  $\alpha: F \rightarrow G$ is an order-preserving  morphism. The following are equivalent:

(1)  $\alpha$ is a frame morphism.

(2)  For each $u \in \mathcal{O}(X)$, $\alpha_{u}: F(u) \rightarrow
G(u)$ is a frame homomorphism and the following square commutes for any $v
\leq u$ in $\mathcal{O}(X)$
 \[
\xymatrix{ F(u)\ar[r]^{\alpha_{u}}&  G(u)\\
           F(v)\ar[r]^{\alpha_{v}}\ar[u]^{f_{uv}}
           & G(v)\ar[u]^{g_{uv}}}
\] where $f_{uv}: F(v) \rightarrow F(u)$ and $g_{uv}: G(v) \rightarrow G(u)$ are the left adjoint of the restriction maps $F(u) \rightarrow
F(v)$ and $G(u) \rightarrow G(v)$ respectively. \end{lemma}

We write $FrmSh(X)$ for the category of all frame sheaves on $X$ and frame morphisms. We now show that the frame sheaves category $FrmSh(X)$ is equivalent to the  category  $\mathcal{O}(X)/ Frm$ of frames under $\mathcal{\mathcal{O}}(X)$, this shows that the frame sheaves on $X$ are just the internal frames in the localic topos $Sh(X)$.

Let $f: \mathcal{O}(X) \rightarrow L$ be a fame homomorphism. We define a sheaf $\Phi(F)$ on $X$ such that $ \Phi(f)(u) = \{x \in L\mid x \leq f(u)\}$, and for $v \leq u$ in $\mathcal{O}(X)$, the restriction map $\Phi(F)(u) \rightarrow \Phi(F)(v)$ sends an element $x \in \Phi(F)(u)$ to an element $x \wedge f(v) \in \Phi(F)(v)$. It is clear that $\Phi(f)$ is a fame sheaf. Suppose we are given a commutative diagram of frame homomorphisms
\[
\xymatrix{&\mathcal{O}(X)\ar[dl]_{f}\ar[dr]^{g}&\\
L\ar[rr]^{h}&&M}
\]
Then for any $u \in \mathcal{O}(X)$, $h$ can be restricted to a frame homomorphism $\Phi(h)(u): \Phi(f)(u) \rightarrow \Phi(g)(u)$, and for  $v \leq u$,  the following square commutes \[
\xymatrix{ \Phi(f)(u)\ar[r]^{\Phi(h)(u)}&  \Phi(g)(u)\\
           \Phi(f)(v)\ar[r]^{\Phi(h)(v)}\ar[u]^{}
           & \Phi(g)(v)\ar[u]^{}}
\] Hence $\Phi(f): \Phi(f) \rightarrow\Phi(g)$ is a frame morphism between frame sheaves. This shows that $\Phi: \mathcal{O}(X) /Frm \rightarrow FrmSh(X)$ is a functor.

\begin{theorem} The frame sheaves category  $FrmSh(X)$  is equivalent to the  category $\mathcal{O}(X)/ Frm$ of frames under $\mathcal{\mathcal{O}}(X)$. \end{theorem}

{\bf Proof}\, Suppose $F$ is a frame sheaf on $X$, write $\Psi(F) = F(1_{X})$ where $1_{X}$ be the largest element of $\mathcal{O}(X)$. For $u \in \mathcal{O}(X)$, let $f: \mathcal{O}(X) \rightarrow \Psi(F), u \mapsto l_{u 1_{X}}(1_{F(u)})$ where $l_{u 1_{X}}: F(u) \rightarrow F(1_{X})$ be the left adjoint of the restriction map $F(1_{X}) \rightarrow F(u)$ and $1_{F(u)}$ be the largest element of $F(u)$. We first show that
$f: \mathcal{O}(X) \rightarrow \Psi(F)$ is a frame homomorphism.

Let $u , v \in \mathcal{O}(X)$.  We have $l_{u\wedge v u}(1_{F(u\wedge v)}) = (l_{v u\vee v}(1_{F(v)}))|_{u}$ by lemma 3.1, it implies that $l_{u\wedge v u}(1_{F(u\wedge v)}) = (l_{v 1_{X}}(1_{F(v)}))|_{u}$. This implies that $l_{u\wedge v 1_{X}}(1_{F(u\wedge v)}) = l_{u 1_{X}}((l_{v 1_{X}}(1_{F(v)}))|_{u}) = l_{u 1_{X}}((l_{v 1_{X}}(1_{F(v)}))|_{u} \wedge 1_{F(u)}) = l_{v 1_{X}}(1_{F(v)}) \wedge l_{u 1_{X}}(1_{F(u)})$ by proposition 3.7. Hence $f$ preserves finite meets. Suppose $\{u_{i} \mid i \in I\} \subseteq \mathcal{O}(X)$, write $t = \bigvee l_{u_{i}\bigvee u_{i}}(1_{F(u_{i})})$. Then $t|_{u_{i}} \geq (l_{u_{i}\bigvee u_{i}}(1_{F(u_{i})})|_{u_{i}} = 1_{F(u_{i})}$. Thus $t = 1_{F(\bigvee u_{i})}$ by sheaf axiom. This implies that $l_{\bigvee u_{i} 1_{X}}(1_{F(\bigvee u_{i})}) = l_{\bigvee u_{i} 1_{X}}(\bigvee l_{u_{i}\bigvee u_{i}}(1_{F(u_{i})})) = \bigvee l_{u_{i}1_{X}}(1_{F(u_{i})}))$ since $l_{\bigvee u_{i} 1_{X}}$ preserves joins.

Suppose $F$ and $G$ are frame sheaves on $X$ and $\alpha: F \rightarrow G$ is a frame morphism. Write $\Psi(\alpha) = \alpha_{1_{X}}: \Psi(F) = F(1_{X}) \rightarrow \Psi(G) = G(1_{X})$. The above argument shows that $\Psi: FrmSh(X)\rightarrow \mathcal{O}(X) /Frm $  is a functor.

We now show that $\Phi \Psi$ is isomorphic to the identity on $ FrmSh(X)$ and $\Psi \Phi$ is isomorphic to the identity on $ \mathcal{O}(X) /Frm $.

 Suppose $F$ is a frame sheaf on $X$. For $v \leq u$ in $\mathcal{O}(X)$,  $l_{vu}: F(v) \rightarrow F(u)$ is the left adjoint
 of restriction map $F(u) \rightarrow F(v)$ . We have $l_{vu}(y \wedge z) = l_{vu}((l_{vu}(y))|_{v} \wedge z) = l_{vu}(y) \wedge l_{vu}(z)$ by proposition 3.7.  Thus
  $l_{vu}: F(v) \rightarrow F(u)$ is a frame mono-homomorphism, in particular, $l_{u 1_{X}}: F(u) \rightarrow F(1_{X})$ is a  frame mono-homomorphism for any $u \in \mathcal{O}(X)$. Moreover, if $x \leq  l_{u 1_{X}}(1_{F(u)})$, then $x = x \wedge  l_{u 1_{X}}(1_{F(u)}) =  l_{u 1_{X}}(x|_{u} \wedge 1_{F(u)}) = l_{u 1_{X}}(x|_{u})$. This shows that $l_{u 1_{X}}: F(u) \rightarrow \Phi \Psi(F)(u) = \downarrow l_{u 1_{X}}(1_{F(u)})$ is an isomorphism of frames. For naturality of the isomorphism,
  we note that $l_{u 1_{X}}(x) \wedge l_{v 1_{X}}(1_{F(v)}) = l_{v 1_{X}}(l_{u 1_{X}}(x)|_{v} \wedge 1_{F(u)}) = l_{v 1_{X}}(x|_{v})$. Hence $\Phi \Psi$ is isomorphic to the identity on $ FrmSh(X)$.  It is clear that $\Psi \Phi $ is isomorphic to the identity on $ \mathcal{O}(X) /Frm $.  $\hfill\Box$

\section{Partially Ordered Sheaf Locales}

It is well known that the category $Sh(X)$ of sheaves on a locale
$X$ is equivalent to the slice category $LH/X$ of locales and
local homeomorphisms over $X$ (see Johnstone [5]).  In this
section, we first give an explicit description of the construction
of the associated sheaf locales and show directly  that the
category $Sh(X)$ of sheaves on a locale $X$ is equivalent to the
slice category $LH/X$ of locales and local homeomorphisms over
$X$, then we give  characterizations of partially ordered sheaf
locales and complete partially ordered sheaf locales respectively.

Recall a localic map $f: X \rightarrow Y$ is said to be local
homeomorphism if $X$ can be covered by open sublocales $U$ for which
the composite $U \rightarrowtail X \rightarrow Y$ is isomorphic to
the inclusion of an open sublocale of $Y$. We will write $LH$ for
the category of locales and local homeomorphisms.

 Let $X$
be a locale, $P \in [\mathcal{O}(X)^{op}, Set]$ be a presheaf on
$X$.  For $s_{i}\in P(u_{i}), u_{i}\in \mathcal{O}(X),
i=1,\cdot\cdot\cdot, n$, write $\epsilon_{P}(s_{1}, \cdot\cdot\cdot,
s_{n})= \bigvee\{u\leq u_{1} \wedge \cdot\cdot\cdot \wedge u_{n}\mid
{s_{i}}\mid_{u}= {s_{j}}\mid_{u}, 1\leq i,j\leq n\}$. If $P$ is a
sheaf, then it is clear that $\epsilon_{P}(s_{1}, \cdot\cdot\cdot,
s_{n})$ is the largest open sublocale $u$ such that ${s_{i}}\mid_{u}
= {s_{j}}\mid_{u}$ for every $i, j = 1, \cdot \cdot \cdot, n$.

\begin{lemma} (1)\, $\epsilon_{P}(s_{1}, \cdot\cdot\cdot,
s_{n}, t_{1}, \cdot\cdot\cdot, t_{m}) \leq
\epsilon_{P}(s_{1},\cdot\cdot\cdot, s_{n}) \wedge
\epsilon_{P}(t_{1},\cdot\cdot\cdot, t_{m}) $.

(2)\, $\epsilon_{P}(s_{1}, \cdot\cdot\cdot, s_{n}, t_{1},
\cdot\cdot\cdot, t_{m}) = \epsilon_{P}(s_{1},\cdot\cdot\cdot, s_{n})
\wedge \epsilon_{P}(t_{1},\cdot\cdot\cdot, t_{m}) $ for any
$s_{i}\in P(u_{i}), t_{j} \in P(v_{j})$ where $u_{i}, v_{j} \in
\mathcal{O}(X)$ if and only if $P$ is the terminal object of
$[\mathcal{O}(X)^{op}, Set]$ .\end{lemma}

Let $X$ be a locale, $P \in [\mathcal{O}(X)^{op}, Set]$.
 We define $\Lambda(P)$ the frame of all functions $f: \coprod
P(u)\rightarrow \mathcal{O}(X)$ with pointwise partial order such
that $\forall f \in \Lambda(P)$ satisfying:

$(1_{\Lambda})\,\, f(s)\leq u$ for $s\in P(u)$;

$(2_{\Lambda})\,\,  f(s) \wedge \epsilon_{P}(s, t) = f(t) \wedge
\epsilon_{P}(s, t)$

 where $ \coprod P(u)$ be the disjoint
union of all $P(u)$ for $u \in \mathcal{O}(X)$. Equivalently,
$\Lambda(P)$ be the subframe of the frame product $\prod_{s\in
\coprod P(u)}\downarrow u_{s}$, where $u_{s} = u$ for $s\in P(u)$,
such that each element $(x_{s})$ of $\Lambda(P)$ satisfying $x_{s}
\wedge \epsilon_{P}(s, t) = x_{t} \wedge \epsilon_{P}(s,t)$ for all
$s, t\in \coprod P(u)$. Let
\[p^{\ast}: \mathcal{O}(X)\rightarrow \Lambda(P),\,\, x\mapsto (x \wedge u_{s}), s\in P(u_{s}), u_{s}\in \mathcal{O}(X)\]
then it is clear that $p^{\ast}$ is a frame homomorphism.

\begin{lemma}  $p: \Lambda(P) \rightarrow X$ is a local homeomorphism.
\end{lemma}

{\Proof}\, Suppose $s\in P(u)$ and $ p_{s}^{\ast}: \Lambda(F)
\rightarrow  \downarrow u$ be the s'th projection. It is clear that
$p_{s}^{\ast}$ is surjective since $(x \wedge u_{s}) \in \Lambda(F)$
for any $x \leq u$. Given $(x_{t}), (x_{t}') \in \Lambda(P)$, we
have $x_{s}'\leq x_{s} \Leftrightarrow x_{t}' \wedge \epsilon_{P}(s,
t) = x_{s}' \wedge \epsilon_{P}(s, t) \leq x_{s} \wedge
\epsilon_{P}(s, t) = x_{t} \wedge \epsilon_{P}(s, t) \leq x_{t}$ for
all $t\in \coprod P(u) \Leftrightarrow x_{t}' \leq \epsilon_{P}(s,
t)\rightarrow x_{t}$ for all $t  \in \coprod P(u) \Leftrightarrow
(x_{t}') \leq (\epsilon_{P}(s, t)) \rightarrow
 (x_{t})$ . Hence the nucleus induced by $p_{s}^{\ast}$ is just the
open nucleus $(\epsilon_{P}(s, t))\rightarrow (\, )$, that is
$\downarrow u$ is isomorphic to the open sublocale
$(\epsilon_{P}(s,t))$ of $\Lambda(P)$. This shows that the set of
opens $\{(\epsilon_{P}(s,t))\mid s \in \coprod P(u)\}$ of
$\Lambda(P)$ form a cover of $\Lambda(P)$ such that each composite
$(\epsilon_{P}(s,t)) \rightarrowtail \Lambda(P) \rightarrow X$ is
isomorphic to an open inclusion. $\hfill\Box$

\begin{proposition} $\Lambda: [\mathcal{O}(X)^{op}, Set] \rightarrow LH/X$ is a
functor.\end{proposition}

{\Proof}\, Let $F, G\in [\mathcal{O}(X)^{op}, Set]$ and $\alpha:
F\rightarrow G$ be a natural transformation. define
$\Lambda(\alpha)^{\ast}: \Lambda(G) \rightarrow \Lambda(F)$ as
\[(x_{t}) \mapsto (x_{s}'),   x_{s}' = x_{t}\,\, for\,\, s\in F(u), t\in
G(u)\,\, and\,\, \alpha_{u}(s) = t\] Note that $\epsilon_{F}(s_{1},
s_{2}) \leq \epsilon_{G}(t_{1}, t_{2})$ for all $s_{1}\in F(u),
s_{2}\in F(v), t_{1}\in G(u), t_{2}\in G(v)$ with $t_{1} =
\alpha_{u}(s_{1}), t_{2} = \alpha_{v}(s_{2})$ by the natural
transformation of $\alpha$. So $x_{s_{1}}' \wedge
\epsilon_{F}(s_{1}, s_{2}) = x_{s_{2}}' \wedge \epsilon_{F}(s_{1},
s_{2})$ for all $s_{1}, s_{2}\in \coprod F(u)$. This shows that
$\Lambda(\alpha)^{\ast}$ is well defined. It is clear that
$\Lambda(\alpha)^{\ast}$ is a frame homomorphism such that $
\Lambda(\alpha)^{\ast} g^{\ast} = f^{\ast}$, i.e. the following
diagram commutes:
\[
\xymatrix{
 \Lambda(F) \ar[rr]^{\Lambda(\alpha)}\ar[dr]_<<<<<{f}& & \Lambda(G)\ar[dl]^<<<<<{g}\\
 &  X  &}
\]

If $\alpha: F\rightarrow G, \beta: G\rightarrow H$ are natural
transformations, then it is clear $\Lambda(\beta \alpha) =
\Lambda(\beta) \Lambda(\alpha)$. Hence $\Lambda:
[\mathcal{O}(X)^{op}, Set] \rightarrow LH/X$ is a functor.
$\hfill\Box$

 Let $X$ be a locale. Recall the cross-sections functor
$\Gamma: Loc/ X\rightarrow Sh(X)$ defined in [5]: given a locale $p:
E\rightarrow X$ over $X$, $\Gamma(p)$ is the sheaf such that
$\Gamma(p)(u)$ be the set of all continuous sections of $p$ over
$u$, where $u\in \mathcal{O}(X)$, i.e. localic maps $s: \downarrow u
\rightarrow E$ such that the composite $p s$ is the inclusion
$\downarrow u \rightarrowtail X$.

\begin{proposition} For any locale $X$, $\Lambda:
[\mathcal{O}(X)^{op}, Set] \rightarrow Loc/ X$ is left adjoint to
the functor $\Gamma: Loc/ X \rightarrow [\mathcal{O}(X)^{op}, Set]$.
The unit $\eta: 1_{[\mathcal{O}(X)^{op}, Set]} \rightarrow \Gamma
\Lambda$ is defined as:
\[\eta_{P_{u}}: P(u)\rightarrow \Gamma \Lambda(P)(u),\,\, s \mapsto
p_{s}\] while the counit $\varepsilon: \Lambda \Gamma \rightarrow
1_{Loc/ X}$ is defined such that for any locale $f: Y\rightarrow X$
over $X$,
\[\varepsilon_{f}^{\ast}: Y \rightarrow \Lambda \Gamma(f), \,\, y \mapsto (s^{\ast}(y))\]
\end{proposition}

{\Proof}\, We first show that $\eta$ and $\varepsilon$ are both
natural transformations.

Suppose $P \in [\mathcal{O}(X)^{op}, Set]$. Note that for $v \leq u$
in $\mathcal{O}(X)$ and $s \in P(u)$, we have $\epsilon_{P}(s,
s\mid_{v}) = v$. Hence $p_{s} i = p_{s\mid_{v}}$ where $i: v
\rightarrowtail u$ is the inclusion. This shows $\eta_{P}: P
\rightarrow \Gamma \Lambda (P)$ is a natural transformation. It is
readily to verify that $\eta_{P}: P \rightarrow \Gamma \Lambda (P)$
is natural for $P$.

For natural transformation $\varepsilon: \Lambda \Gamma \rightarrow
1_{Loc/ X}$ it is clear.

Next we observe that $\eta$ and $\varepsilon$ have the property that
both composites
 \[\Gamma \stackrel{\eta \Gamma}{\rightarrow} \Gamma \Lambda \Gamma
 \stackrel{\Gamma \varepsilon}{\rightarrow} \Gamma,\,\, \Lambda \stackrel{\Lambda \eta}
 {\rightarrow} \Lambda \Gamma \Lambda \stackrel{\varepsilon \Lambda}{\rightarrow} \Lambda\]
are identities. For the first composite, given a locale $f:
Y\rightarrow X$ over $X$ and  $s\in \Gamma(f)(u)$,
$\eta_{\Gamma(f)_{u}}$ sends $s$ to the s'th projection $p_{s}: u
\rightarrow \Lambda \Gamma(f)$, and then
$\Gamma(\varepsilon_{f})_{u}$ sends it to $ \varepsilon_{f} p_{s} =
s$. Similarly, for a presheaf $P$ on $X$,  the second composite
first sends $y = (x_{t})\in \Lambda(P)$ to $(s^{\ast}(y))\in \Lambda
\Gamma \Lambda (P)$ by $\varepsilon_{\Lambda(P)}^{\ast}$ and then
sends it to $(p_{t}^{\ast}(y)) = (x_{t})$ by $\Lambda
(\eta_{P})^{\ast}$. $\hfill\Box$

\begin{corollary} The functor $\Lambda$ and $\Gamma$ in Proposition 5.2
restrict to an equivalence of categories
\[Sh(X) \rightleftarrows LH/ X\]
Moreover, $Sh(X)$ is a reflective subcategory of
$[\mathcal{O}(X)^{op}, Set]$, and $LH/ X$ is a coreflective
subcategory of $Loc/ X$.
\end{corollary}

{\Proof}\, Let $f: Y\rightarrow X$ be a local homeomorphism. Then we
have a cover $\bigvee y_{i} = Y$ such that each composite
$\downarrow y_{i} \stackrel{s_{i}}{\rightarrowtail} Y
\stackrel{f}{\rightarrow} X$ is an open inclusion. For $y, y'\in
\mathcal{O}(Y)$ with $\varepsilon_{f}^{\ast}(y) =
\varepsilon_{f}^{\ast}(y')$, we have $y = \bigvee y \wedge y_{i} =
\bigvee s_{i}^{\ast}(y) = \bigvee s_{i}^{\ast}(y') = y'$. Hence
$\varepsilon_{f}^{\ast}$ is one to one. For each $(x_{s}) \in
\Lambda \Gamma (f)$, there exist $y_{s}\in \mathcal{O}(Y)$ such that
$x_{s} = s^{\ast}(y_{s})$ for every $s\in \Gamma(f)$. But
$s^{\ast}(y_{s}) \wedge \epsilon(s, t) = t^{\ast}(y_{t}) \wedge
\epsilon(s, t) = s|_{\epsilon(s ,t)}^{\ast}(y)$ for $s, t\in
\Gamma(f)$ implies that $y_{s} = y_{t} = y$. Thus $(x_{s}) =
(s^{\ast}(y))$ for some $y\in \mathcal{O}(Y)$. This shows that
$\varepsilon_{f}^{\ast}: Y \rightarrow \Lambda \Gamma(f)$ is an
isomorphism.

Suppose we are given a sheaf $F$ on $X$, $u$ is an open of $X$ and
$\lambda: u \rightarrow \Lambda(F)$ be a continuous section over $u$.
Then for any $s \in F(v)$ and $s' \in F(v')$, we have
$\lambda^{\ast}((\epsilon_{F}(s, t))) \wedge
\lambda^{\ast}((\epsilon_{F}(s', t))) =
\lambda^{\ast}(((\epsilon_{F}(s, t))\wedge (\epsilon_{F}(s', t))))
\leq \lambda^{\ast}((\epsilon_{F}(s, s')) \leq
\lambda^{\ast}((\epsilon_{F}(s, s'))\wedge u_{s}) = \epsilon_{F}(s,
s')$, and also $\bigvee_{s\in \bigcup F(v)}
\lambda^{\ast}((\epsilon_{F}(s, t))) = 1_{u}$. Hence
$\{s|_{\lambda^{\ast}((\epsilon_{F}(s, t))}\mid s\in \bigcup
F(v)\}$ form a compatible family, and so they patch to an unique
element $\bar{s}\in F(u)$ such that
$\bar{s}|_{\lambda^{\ast}((\epsilon_{F}(s, t))} =
s|_{\lambda^{\ast}((\epsilon_{F}(s, t))}$ for any $s\in \bigcup
F(v)$. By the pullback property of the following diagram:
\[
\xymatrix{\downarrow \lambda^{\ast}((\epsilon_{F}(s,
t)))\ar[r]^{}\ar[d]_{}& \downarrow (\epsilon_{F}(s, t))
\ar[d]^{}\\
      \downarrow u\ar[r]^{\lambda}& \Lambda(F)  }
\]
we know that $\lambda$ must be the $\bar{s}$'th projection
$p_{\bar{s}}$. Hence $\eta_{P_{u}}: P(u)\rightarrow \Gamma
\Lambda(P)(u)$ is a bijection. $\hfill\Box$

Let $f: Y \rightarrow X$ be a locale morphism and $F$ be the corresponding frame sheaf under the equivalence of theorem 3.1. Recall that $F$ is called spatial in the localic topos $Sh(X)$  if  $F$ is isomorphic to a subframe of a power object, equivalently, if there exists a locale epimorphism $h: E \rightarrow Y$ such that the composition $f h$ is a local homeomorphism. By proposition 5.2, we know that $F$ is spatial if and only if the counit $\varepsilon: \Lambda \Gamma \rightarrow
1_{Loc/ X}$ is an epimorphism.  Thus we have the following result.

\begin{corollary} Let $f: Y \rightarrow X$ be a locale morphism.  $f: Y \rightarrow X$ corresponds to a spatial frame if and only if for $\forall y, z \in \mathcal{O}(Y)$, $y \not= z$, there exists a continuous section $s$ over some open $u \in \mathcal{O}(X)$ such that $s^{\ast}(y) \not = s^{\ast}(z)$. \end{corollary}

Let  $f: Y \rightarrow X$ be a local homeomorphism and $u , v \in
\mathcal{O}(X)$ be two opens with $v \leq u$. If $s \in
\Gamma(f)(u)$ be a section over $u$, then we have section $s|_{v}$ over
$v$ such that for any $y \in \mathcal{O}(Y)$, $s|_{v}^{\ast}(y) =
s^{\ast}(y) \wedge v$. We call it the restriction of $s$ to $v$. Now
we define the concept of partially ordered sheaf locales in such a
way as to make the equivalence $Sh(X) \rightleftarrows LH/ X$ in
corollary 5.1 true for the category of partially ordered sheaves.

\begin{definition} Let $f: Y \rightarrow X$ be a local
homeomorphism. $f$ is said to be a partially ordered sheaf locale if

(POSL1) For any $u \in \mathcal{O}(X)$, the set $\Gamma(f)(u)$ of
all sections on $u$ is a partially ordered set.

(POSL2) If $v \leq u$ in $\mathcal{O}(X)$,  $s , t \in \Gamma(f)(u)$
be two section on $u$ with $s \leq t$, then $s|_v \leq t|_{v}$ in
$\Gamma(f)(v)$.

(POSL3) If we have a cover $u = \bigvee u_{i}$ in $\mathcal{O}(X)$,
and two sections $s, t \in \Gamma(f)(u)$ such that $s|_{u_{i}} \leq
t|_{u_{i}}$ in $\Gamma(f)(u_{i})$ for each $u_i$ then  $s \leq t$ in
$\Gamma(f)(u)$. \end{definition}

Given a map of locales over $X$ with $f:Y \rightarrow X$ and $g: Z
\rightarrow X$ are both local  homeomorphisms
\[
\xymatrix{
 X \ar[rr]^{\phi}\ar[dr]_<<<<<{f}& & Y \ar[dl]^<<<<<{g}\\
 &  Z  &}
\] We call $\phi$ an order-preserving map if for each $u \in
\mathcal{O}(X)$, the map \[\phi_{u}: \Gamma(f)(u) \rightarrow
\Gamma(g)(u),  s \mapsto \phi s\] is an order-preserving map between
posets. We write $POLH/ X$ for the category of all partially ordered
sheaf locales over $X$ and order-preserving maps, and $\mathfrak{POS}_{X}$ for the category
of all posheaves on $X$ and order-preserving morphisms.  Then we have

\begin{corollary} The functors $\Lambda$ and $\Gamma$ in Proposition 4.2
restrict to an equivalence of categories
\[\mathfrak{POS}_{X} \rightleftarrows POLH/ X\].
\end{corollary}

To characterize complete partially ordered sheaves in terms of
sheaf locales, we need to introduce the concepts of complete
partially ordered sheaf locales.

\begin{definition} Let $f: Y \rightarrow X$ be a local
homeomorphism. $f$ is said to be a complete partially ordered sheaf
locale if

(CPOSL1) For any $u \in \mathcal{O}(X)$, the set $\Gamma(f)(u)$ of
all sections on $u$ is a complete lattice.

(CPOSL2) If $v \leq u$ in $\mathcal{O}(X)$,  then the restriction
map $\Gamma(f)(u) \rightarrow \Gamma(f)(v), s \mapsto s|_v$ is
surjective,  and preserves arbitrary joins and meets.

(CPOSL3) If we have a cover $u = \bigvee u_{i}$ in $\mathcal{O}(X)$,
and two sections $s, t \in \Gamma(f)(u)$ such that $s|_{u_{i}} \leq
t|_{u_{i}}$ in $\Gamma(f)(u_{i})$ for each $u_i$ then  $s \leq t$ in
$\Gamma(f)(u)$. \end{definition}

We write $CPOLH/ X$ for the category of all complete partially
ordered sheaf locales over $X$ and order-preserving maps,  $\mathfrak{CPOS}_{X}$ for the category of all complete
posheaves on $X$ and order-preserving morphisms. Then we have

\begin{corollary} The functors $\Lambda$ and $\Gamma$ in Proposition 4.2
restrict to an equivalences of categories \[\mathfrak{CPOS}_{X}
\rightleftarrows CPOLH/ X\]
\end{corollary}